\documentclass[12pt]{article}
\usepackage[margin=1in]{geometry}

\usepackage[utf8]{inputenc}
\usepackage{amsmath}
\usepackage{hyperref}
\usepackage{amssymb}
\usepackage{mathtools}
\usepackage{graphicx}
\newcommand{\bfnu}{\boldsymbol{\nu}}

\newcommand{\bfxi}{\boldsymbol{\xi}}
\newcommand{\bfmu}{\boldsymbol{\mu}}
\newcommand{\bfphi}{\boldsymbol{\phi}}
\newcommand{\bfw}{\mathbf{w}}
\newcommand{\bfE}{\mathbf{E}}
\newcommand{\bfH}{\mathbf{H}}
\newcommand{\bfU}{\mathbf{U}}
\newcommand{\bfu}{\mathbf{u}}
\newcommand{\bfd}{\mathbf{d}}
\newcommand{\bfx}{\mathbf{x}}
\newcommand{\bfg}{\mathbf{g}}
\newcommand{\0}{{\mathbf{0}}}
\newcommand{\bft}{\mathbf{t}}
\newcommand{\bfp}{\mathbf{p}}
\newcommand{\bfq}{\mathbf{q}}
\newcommand{\bbS}{\mathbb{S}}
\newcommand{\cS}{\mathcal{S}}
\newcommand{\bfv}{\mathbf{v}}

\newcommand{\DE}{{\mathbb R}^3\setminus \overline{D}}

\def\hd1{{\H^{1/2}(\div_{\partial B},\partial B)}}

\newcommand{\sL}{{L}^2_t}
\newcommand{\hbfx}{\hat{\mathbf{x}}}
   
\DeclareMathOperator{\curl}{curl}
\DeclareMathOperator{\bfcurl}{\mathbf{curl}}
\DeclareMathOperator{\Div}{div}
\usepackage{color}

\newtheorem{theorem}{Theorem}
\newtheorem{lemma}{Lemma}
\newtheorem{definition}{Definition}
\newtheorem{assumption}{Assumption}
\newtheorem{prop}{Proposition}

\title{Target Signatures for  Anisotropic Screens in Electromagnetic Scattering}

\author{Fioralba Cakoni\thanks{Department of Mathematics, Rutgers University, New Brunswick, NJ 08903, USA.
    ({\tt fc292@math.rutgers.edu})}
\and Peter Monk\thanks{Department of Mathematical Sciences, University of Delaware, Newark,  Delaware, USA. ({\tt monk@udel.edu})}}

\begin{document}

\maketitle

 \begin{abstract}
Anisotropic thin sheets of materials possess intriguing properties because of their ability to modify the phase, amplitude and polarization of incident waves.  Such sheets are usually modeled by imposing transmission conditions of resistive or conductive type on a surface called a screen.  We start by analyzing this model, and show that the standard passivity conditions can be slightly strengthened to provide conditions under which the forward scattering problem has a unique solution.  We then turn to the inverse problem and suggest a target signature for monitoring such films.  The target signature is based on  a modified far field equation obtained by subtracting an artificial far field operator for scattering by a closed surface containing the thin sheet and parametrized by an artificial impedance.  We show that this impedance gives rise to an interior eigenvalue problem, and these eigenvalues can be determined from the far field pattern, so functioning as target signatures.  We prove uniqueness for the inverse problem, and give preliminary numerical examples illustrating our theory.
 \end{abstract}

\section{Introduction}
Ultra-thin sheets of materials such as graphene have been the subject of intensive research for several decades~\cite{su141912336} because they can be tuned  to modify the phase, amplitude and polarization of incident waves.  More recently, the possibility of using thin sheets of meta-materials has expanded the 
range of possible behaviors of the sheet to include anisotropic surface  surface properties (see for example \cite{Alu,alu16,Alu2,Alu3,7565471}).  Such ultra-thin structures, hereafter called screens,  are usually modeled by imposing transmission conditions across the screen using a suitable optical conductivity tensor~\cite{alu16}.  This model can be derived as a limiting case of a thin penetrable material layer~\cite{DHJ_12,irene_16} as the thickness tends to zero. The resulting transmission problem contrasts to models of thin materials that have prescribed boundary conditions (for example \cite{CEA2015,ola23}), so that new theory needs to be derived.

\noindent
The first step in this paper is to study a general model for forward scattering by ultra-thin screens. More precisely, assuming a complete description of the screen, we want to predict how it scatters incoming radiation.  We  prove that the forward problem is well posed in the important case of a uniaxial passive metasurface, so connecting a strengthened form of the usual assumptions of passivity~\cite{alu16} to coercivity of certain sesquilinear forms, and hence using Fredholm theory, to the existence of a unique solution to the forward problem.  We then move on to the inverse problem of detecting changes in the material properties of the isotropic or anisotropic screens using target signatures.  In this context, target signatures are discrete quantities that can be computed from scattering data.  Changes
in these quantities could then be used to monitor or detect changes in the screen.  Typically these quantities are eigenvalues of an interior problem.  They arise by modifying the far field operator using an auxiliary far field operator generated by a suitable parameter dependent 
problem. Building on previous work for electromagnetism in two dimensions~\cite{Screens,heejin22}, we suggest a new target signature derived by
considering the injectivity of a modified far field operator for the 3D Maxwell problem. We characterize the target signatures as eigenvalues of an interior problem where we suppose that the screen covers a part of the boundary of an artificial closed bounded domain in $\mathbb{R}^3$ on which the eigenvalue problem is defined.  This target signature is simpler than our previous 2D signatures for thin screens in that the auxiliary scattering problem that contributes to the modified far field operator is independent of the details of the conducting screen.

\noindent
The paper is structured as follows.  In Section~\ref{model} we introduce the function spaces used on this paper, and present the forward problem of scattering by a known screen.  We derive an existence theory for such problems that encompasses models reported in the literature (e.g.~\cite{alu16}).   In Section~\ref{inverse} we discuss the inverse problem of determining the surface impedance from far field data, and prove a uniqueness theorem for the problem suggesting that the data we use for target signature is rich enough to characterize the screen.  We then define the modified far field operator and the target signatures for this paper. We prove a relationship between the target signatures and injectivity of the modified far field operator.  In  Section~\ref{sigma-stek} we study the eigenvalue problem related to our target signatures called the $\Sigma$-Steklov eigenvalue problem. Section~\ref{numerics}  presents a discussion on  the determination of $\Sigma$-Steklov eigenvalues from far field data, and shows  some preliminary numerical results illustrating our theory.

\section{Notation and the Forward Problem}\label{model}
We start this section by summarizing the function spaces needed for this paper.  Then we move on to discuss the forward scattering problem for a thin resistive or conductive screen.  This problem will underly our discussion of the inverse problem.

\noindent
The thin screen occupies a region $\Gamma\subset {\mathbb R}^3$ denoting a piecewise smooth, compact, open  two dimensional  manifold with boundary. We assume that $\Gamma$ is simply connected and non self-intersecting such that  it can be embedded as part of a piece-wise smooth closed boundary $\partial D$ circumscribing a bounded 
connected region $D\subset {\mathbb R}^3$ having connected complement. This determines two sides of $\Gamma$ and  we choose the positive side using the unit normal vector $\bfnu$ on $\Gamma$ that coincides with the normal direction outward of $D$. To be able to precisely define the scattering problem and for later use we recall the definition of several Sobolev spaces:

\noindent
\subsection{Function spaces}
Let ${\cal Y}$ be a domain in $\mathbb{R}^3$ then recall the standard space of curl conforming vector functions on ${\cal Y}$ 
$$
H(\curl, {\cal Y}):= \left \{{\bfu}\in
(L^2({\cal Y}))^3\,:\,\curl{\bfu}\in
(L^2({\cal Y})^3\right \} $$
 and denote by $
H_{loc}(\curl, {\mathbb R}^3)$ the space of $\bfu\in H(\curl,
B_R)$ for all $B_R$ where $B_R$ is a ball centered at the origin with radius $R$ containing $\Gamma$ containing $\Gamma$. Then, using the space of $L^2$ tangential vector fields on $\Gamma$ denoted by $\sL(\Gamma)$, we define the Sobolev space
$$
X(\curl, B_R) :=\{\bfu \in H(\curl, B_R)\,:\,
\bfu_T \in \sL(\Gamma)\},
$$
endowed with the natural norm
$$
\|\bfu\|^2_{X(\curl, B_R)}:=\|\bfu\|_{H(\curl,B_R)}^2+\|\bfu_T \|^2_{L^2(\Gamma)}
$$
where $\bfu_T=(\bfnu\times\bfu)\times\bfnu$. Next let $D$ be a bounded region in ${\mathbb R}^3$ with piecewise smooth boundary $\partial D$ such that $\Gamma\subset \partial D$, chosen such that the positive side of $\Gamma$ coincide with the outward direction on $\partial D$. We can also define 
corresponding space $H_{loc}(\curl, \DE)$.  Obviously we  also  have
$$
X(\curl, D) :=\{\bfu \in H(\curl, D)\,:\,
\bfu_T \in L_t^2(\Gamma)\},
$$
$$
X(\curl, B_R\setminus {\overline{D}}) :=\{\bfu \in H(\curl, B_R\setminus {\overline{D}})\,:\,
\bfu_T \in L_t^2(\Gamma)\},
$$
and the correspondingly $X_{loc}(\curl, \DE)$. For later use we define additional Sobolev  spaces on the piece-wise smooth  boundary $\partial D$
\[
\begin{array}{ll}
H^s_t(\partial D):=\left\{\bfmu\in H^s(\partial D)^3\,:\,\bfnu\cdot\bfmu=0\hbox{  a.e. on }\partial D \right\}\,, \\[2ex]
H^s(\Div_{\partial D},\partial D ):=\left\{\bfmu\in H^s_t(\partial D)\,:\,\Div_{\partial D} \bfmu\in H^s(\partial D)\right\}\,, \\[2ex]
H^s(\Div_{\partial D}^0,\partial D):=\left\{\bfmu\in H^s(\Div_{\partial D},\partial D)\,:\,\Div_{\partial D} \bfmu =0 \hbox{  on }\partial D\right\}\,, \\[2ex]
H^{-1/2}(\curl_{\partial D},\partial D):=\left\{\bfmu\in H^{-1/2}_t(\partial D)\,:\,\curl_{\partial D} \bfmu\in H^{-1/2}(\partial D)\right\}\,,
\end{array}
\]
where $\curl_{\partial D}$ and $\Div_{\partial D}$ are the surface scalar $\curl$ and divergence operator, respectively, and $s\in\mathbb {R}$. In addition we will denote by $\bfcurl_{\partial D}$ the surface vectorial $\curl$. We rename the spaces $H^0_t(\partial D)$ and $H^0(\Div_{\partial D},\partial D)$ by  $\sL(\partial D)$ and $H(\Div_{\partial D},\partial D)$, respectively. The space $H^s_t(\partial D)$ is equipped with the standard norm (see, for instance, \cite{monk}), whereas the spaces $H^s(\Div_{\partial D},\partial D)$ and $H^{-1/2}(\curl_{\partial D},\partial D)$ are endowed with their respective natural norms
\[
 \|\bfmu\|_{H^s(\Div_{\partial D},\partial D)}:= \|\bfmu\|_{s,\partial D}^2\,+\,\|\Div_{\partial D}\bfmu\|_{s,\partial D}^2\qquad \qquad \mbox{and}
 \]
\[
\|\bfmu\|^2_{H^{-1/2}(\curl_{\partial D},\partial D)}\,:=\, \|\bfmu\|_{-1/2,\partial D}^2\,+\,\|\curl_{\partial D}\bfmu\|_{-1/2,\partial D} ^2\,.
 \]
 Note that  integration by parts in $H(\curl, D)$ (or $H(\curl, B_R\setminus {\overline{D}})$) defines a duality between the rotated tangential trace in  $H^{-1/2}(\Div_{\partial D},\partial D)$ and the tangential trace in $H^{-1/2}(\curl_{\partial D},\partial D)$. For more details about the norms and properties of this operators, see for instance \cite{monk} for smooth boundaries and \cite{buffa1, buffa2} for Lipschitz boundaries.

\noindent

\subsection{The forward problem}
 We now  rigorously  describe the forward scattering problem.  We first define the time harmonic incident electric field $e^{-i\omega t}\bfE^i(\bfx)$ at angular frequency $\omega$ to be a plane wave, where the spatially dependent part  $\bfE^i$ satisfies the background Maxwell system in all space and is given by
\begin{equation}\label{plane}
\bfE^{i}(\bfx;\kappa, \bfd,\bfp)=\frac{i}{\kappa}\curl \curl \bfp {\rm w}e^{i\kappa \bfd\cdot\bfx}=i\kappa(\bfd\times\bfp)\times \bfd {\rm e}^{i\kappa \bfd\cdot\bfx}.
\end{equation}
Here the unit vector  $\bfd\in\mathbb{R}^3$, $|\bfd|=1$, is the direction of propagation and $\bfp\in\mathbb{C}^3$ is the polarization. To satisfy the background Maxwell's system, we must have $|\bfd|=1$, $\bfp\not=0$ and $\bfd\cdot\bfp=0$. In addition, $\kappa>0$ is the wave number that is related to  the angular frequency $\omega$ of the radiation  by $\kappa=\omega\sqrt{\epsilon_0 \mu_0}$ where $\epsilon_0$ and $\mu_0$ are electric permittivity and magnetic permeability  of the homogenous background medium (free space). Other incident fields can also be used (for example those due to point sources).

 \noindent
Following \cite{irene_16,jin_volakis,Senior}, the electromagnetic properties of a thin screen with central surface $\Gamma$ are described by a matrix valued function $\Sigma$ defined on $\Gamma$. This is a function of position on the screen, its thickness $\delta$, 
 and the physical properties of the screen such as electric permeability, magnetic permittivity and conductivity. We take it  to be a $3\times 3$ piecewise smooth complex valued matrix function of position on $\Gamma$ in order to model an anisotropic screen.  The  tensor $\Sigma$ maps  a vector 
tangential to $\Gamma$ at a point $\bfx\in\Gamma$ to a vector tangential to $\Gamma$ at the same point $\bfx\in\Gamma$. To be more precise, on a smooth face of the surface $\Gamma$ let $\bfnu(\bfx)$  be the smooth outward  unit normal vector function to $\Gamma$ and let  $\hat \bft_1(\bfx)$ and $\hat \bft_2(\bfx)$ be two perpendicular 
vectors in the tangent plane to $\Gamma$ at the point $\bfx$ such that  $\hat \bft_1, \hat \bft_2, \bfnu$ form a right hand coordinative system with origin at  $\bfx$. Using these coordinates,  the matrix valued function $\Sigma(\bfx)$  is represented by the following dyadic expression 
\begin{equation}\label{imp-exp}
\Sigma(\bfx)=\left(\sigma_{11}(\bfx)\hat \bft_1(\bfx)+\sigma_{12}(\bfx)\hat \bft_2(\bfx)\right)\hat \bft_1(x)+\left(\sigma_{21}(\bfx)\hat \bft_1(\bfx)+\sigma_{22}(\bfx)\hat \bft_2(\bfx)\right)\hat \bft_2(\bfx).
\end{equation}
In general, for dispersive thin screens,  $\Sigma:=\Sigma(\bfx, \omega)$ is frequency dependent, but we omit the $\omega$-dependance since our target signatures use scattering data at a single fixed frequency. Note that,  if $\bfxi(\bfx)=\alpha\hat \bft_1(\bfx)+\beta \hat \bft_2(\bfx)$ for some  $\alpha,\beta\in{\mathbb C}$,  then
 $\Sigma(\bfx) \bfxi(\bfx)$ is the tangential vector given by
 $$\Sigma(\bfx)\bfxi(\bfx)=(\alpha \sigma_{11}(\bfx)+\beta\sigma_{21}(\bfx))\hat \bft_1(\bfx)+(\alpha \sigma_{12}(x)+\beta\sigma_{22}(\bfx))\hat \bft_2(\bfx)$$
 and then
\begin{equation}\label{quad}
\overline{\bfxi(\bfx)}^\top \cdot \Sigma(\bfx)\bfxi(\bfx)=|\alpha|^2\sigma_{11}(\bfx)+\overline{\alpha}\beta\,\sigma_{12}(\bfx)+\overline{\beta}\alpha\,\sigma_{21}(\bfx)+|\beta|^2\sigma_{22}(\bfx).
\end{equation}
Generically,  we assume that in the local coordinate system on $\Gamma$, $\Sigma\in (L^{\infty}(\Gamma))^{2\times 2}$ (unless otherwise indicated) thus 
$$\Sigma: \sL(\Gamma)\to \sL(\Gamma) \qquad \mbox{mapping} \qquad {\bfxi}\mapsto \Sigma\bfxi.$$
\noindent
The screen causes a jump in the tangential component of the magnetic field.  To describe this we 
need some notation: for any sufficiently smooth vector field ${\mathbf W}$ defined in ${\mathbb R}^3\setminus \Gamma$ let 
${\mathbf W}^+={\mathbf W}|_{\mathbb{R}^3\setminus \overline D}$ and ${\mathbf W}^-={\mathbf W}|_D$.  In addition, let 
${\mathbf W}^\pm_T=\bfnu\times({\mathbf W}^\pm\times\bfnu)$ on $\Gamma$ the tangential trace from inside and outside. Now, given the screen $\Gamma$ and associated tensor $\Sigma$, as well as the incident field, the forward scattering problem for the screen is to determine the electric field $\bfE$ such that
\begin{subequations}\label{Forward_problem}
\begin{eqnarray}
\curl\curl \bfE-\kappa^2\bfE={\bf 0}&&\qquad \mbox{ in }\mathbb{R}^3\setminus\Gamma,\label{eq:maxwell}\\
\bfE=\bfE^s+\bfE^i&& \qquad \mbox{ in }\mathbb{R}^3\setminus\Gamma,\label{eq:scat}\\
\bfE^+_T=\bfE^-_T&&\qquad \mbox{ on }\Gamma,\label{trans1}\\
\bfnu\times(\curl\bfE^+-\curl\bfE^-)=i\kappa\Sigma\bfE^+_T && \qquad \mbox{ on }\Gamma,\label{trans2}\\
\lim\limits_{|\bfx|\to \infty}\left(\curl \bfE^s\times \bfx-i\kappa|\bfx|\bfE^s \right)=0.&&\label{eq:silver}
\end{eqnarray}
\end{subequations}
Here $\bfE^s$ denotes the scattered electric field, and (\ref{eq:silver}) is the Silver-M{\"u}ller radiation condition
which holds uniformly in all directions $\hbfx=\bfx/|\bfx|$. Equations (\ref{trans1}) and (\ref{trans2})
model the thin anisotropic conductive/resistive thin screen~\cite{irene_16,jin_volakis,Senior}. 

\noindent
First we need to impose conditions on $\Sigma$ in order to guarantee the uniqueness of solutions of the forward problem  (\ref{eq:maxwell})-(\ref{eq:silver}).  Formally, integrating by parts over a ball $B_R$ of radius $R>0$ centered at the origin with $D\subset B_R$,    we have that
\begin{eqnarray*}
&&\int_{B_R} (\curl \bfE^s\cdot\curl\overline{\bfv}-\kappa^2\bfE^s\cdot\overline{\bfv})\,dV-i\kappa\int_{\Gamma}\Sigma \bfE^s_T\cdot\overline{\bfv}_T\,dA \nonumber\\
&&\qquad \qquad +\int_{\partial B_R}\bfnu\times\curl\bfE^s\cdot \overline{\bfv}\,dA=i\kappa\int_{\Gamma}\Sigma\bfE^i_T\cdot\overline{\bfv}_T\,dA. 
\end{eqnarray*}
Now taking ${\bf v}={\bf E^s}$, and choosing $\bfE^i={\mathbf 0}$ we obtain
\begin{eqnarray*}
&&i\kappa\int_{\partial B_R}(\bfnu\times\overline{\bfE}^s)\cdot {\bfH}^s \,dA=\int_{\partial B_R}\bfnu\times\curl\bfE^s\cdot \overline{\bfE}^s \, dA\\
&&\qquad \qquad= \int_{B_R} (|\curl \bfE^s|^2-\kappa^2|\bfE^s|^2\,dV -i\kappa\int_{\Gamma}\Sigma \bfE^s_T\cdot\overline{\bfE}^s_T\,dA
\end{eqnarray*}
Thus  Rellich's Lemma \cite[Theoem 6.10]{CK2019} implies the uniqueness of any solution of (\ref{eq:maxwell})-(\ref{eq:silver}) provided that 
$$\Re\int_{\partial B_R}(\bfnu\times\overline{\bfE}^s)\cdot {\bfH}^s \,dA = - \Re\int_\Gamma \Sigma \bfE^s \cdot\overline{\bfE^s_T} \,dA\leq 0.$$
To provide explicit conditions on the complex valued surface tensor for which the above equality holds, we impose the  condition
\begin{equation}\label{pos}
\Re\left(\overline{\bfxi(\bfx)}^\top \cdot \Sigma(\bfx) \bfxi(\bfx)\right)\geq 0, \;\;\;\;\;  \mbox{$\forall$  complex fields $\bfxi$  tangential to $\Gamma$ a.a. $\bfx\in \Gamma$}
\end{equation}
where the quadratic form is given by (\ref{quad}). Setting 
$$A:=|\alpha|^2\Re(\sigma_{11}), \qquad C:=|\beta|^2\Re(\sigma_{22}), \qquad 2B:=\overline{\alpha}\beta\,(\sigma_{12}+\overline{\sigma}_{21})$$
we see that (\ref{pos}) is satisfied if the Hermitian matrix  $\left(\begin{array}{rrcll} A & B \\ 
\overline{B} & C \end{array}\right)$ is non-negative, i.e. its eigenvalues are non-negative, which is the case provided
\begin{equation}\label{unique}
\Re(\sigma_{11})\geq 0 \qquad \Re(\sigma_{22})\geq 0 \qquad \mbox{and} \qquad \Re(\sigma_{11})\Re(\sigma_{22})\geq 1/4 |\sigma_{12}+\overline{\sigma}_{21}|^2.
\end{equation}
It is easy to see  that (\ref{unique}) can be equivalently written in the following form
$$
\Re(\sigma_{11})\geq 0 \qquad  \Re(\sigma_{22})\geq 0 \qquad \mbox{and} \qquad \Re(\sigma_{11})+\Re(\sigma_{22})\geq |\sigma_{12}+\overline{\sigma}_{21}|
$$
which is customarily  found in  the literature on meta-surfaces \cite{bilow,Alu}.

\noindent
The proof of the existence of the solution  of (\ref{eq:maxwell})-(\ref{eq:silver}) follows the standard approach of \cite{eric,monk}.
Given $\bfE^i$ it is natural to look for the solution $\bfE^s$ of  (\ref{eq:maxwell})-(\ref{eq:silver}) in  $X_{loc}(\curl, B_R)$ (since the tangential component of $\bfE^s$ is continuous across $\Gamma$). Using the exterior Calderon operator, we can reduce the problem to the bounded domain $B_R$.  Then we seek $\bfE^s\in X(\curl, B_R)$ such that 
\begin{eqnarray*}
&&\int_{B_R} (\curl \bfE^s\cdot\curl\overline{\bfv}-\kappa^2\bfE^s\cdot\overline{\bfv})\,dV-i\kappa\int_{\Gamma}\Sigma \bfE^s_T\cdot\overline{\bfv}_T\,dA +i\kappa\int_{\partial B_R}G_e(\hat\bfx\times\bfE^s)\cdot \overline{\bfv}_T\,dA \nonumber\\
&&\qquad  =\int_{\Gamma}i\kappa\eta\bfE^i_T\cdot\overline{\bfv}_T\,dA-i\kappa\int_{\partial B_R}G_e(\hat\bfx\times\bfE^i)\cdot \overline{\bfv}_T\,dA \qquad \forall \;  \bfv\in X(\curl, B_R). \label{varE}
\end{eqnarray*}
 Here $G_e$ is the exterior Calderon operator (c.f. \cite{monk}) which  maps a tangential vector field ${\boldsymbol{\tau}}$ on $\partial B_R$ to $(1/i\kappa)\hat \bfx\times \curl \bfE|_{\partial B_R}$  where the outgoing field $\bfE$ (i.e. satisfying  (\ref{eq:silver})) is a solution of 
$$\nabla \times \bfE-\kappa^2 \bfE= 0\quad \mbox{in}\quad {\mathbb R}^3\setminus \overline{B}_R, \qquad  \hat x \times \bfE = {\boldsymbol{\tau}} \quad \mbox{on}\quad {\partial B}_R.$$ 
The analysis of the terms containing $G_e$ follows exactly the lines of \cite[Theorem 2.3]{ccu} (see also \cite[Theorem 10.2]{monk}) based on a Helmholtz decomposition and  on the fact that the operator $i\kappa G_e$ can be split into a compact part $i\kappa G_e^1$ and a nonnegative part $i\kappa G_e^2$. To avoid repetition, we highlight here the only difference coming from the more general choice of the surface tensor $\Sigma$, which amounts to conditions on $\Sigma$ for which
\begin{eqnarray*}
a({\mathbf W}, {\mathbf W})&=& \int_{B_R}\left(|\curl {\mathbf W}|^2+ |{\mathbf W}|^2\right)\,dA+\kappa \int_{\Gamma}\Im\left(\Sigma {\mathbf W}_T\cdot\overline{{\mathbf W}}_T\right)\,dA \\
&-&i\kappa \int_{\Gamma}\Re\left(\Sigma {\mathbf W}_T\cdot\overline{{\mathbf W}}_T\right)\,dA
\end{eqnarray*}
is coercive in ${X(\curl, B_R)}$, where we have ignored $i\kappa\int_{\partial B_R}G^2_e(\hat\bfx\times\mathbf{W})\cdot \overline{\mathbf{W}}_T\,dA>0$. It is sufficient to find $\theta$ such that, for some $C>0$,
$$\Re\left(e^{i\theta}a({\mathbf W}, {\mathbf W})\right)\geq C\left(\|{\mathbf W}\|_{H({\rm curl},B_R\setminus\overline{\Gamma})}^2+\|{\mathbf W}_T\|^2_{L^2(\Gamma)}\right)$$ which, given (\ref{pos}), is satisfied if for some $0\leq \theta\leq \pi/2$ and $\gamma>0$ constant  and for almost all $\bfx\in \Gamma$,
\begin{equation*}\label{conco}
(\cos \theta) \Re\left(\overline{\bfxi(\bfx)}^\top\cdot \Sigma(\bfx)\bfxi(\bfx)\right) +(\sin \theta) \Im\left(\overline{\bfxi(\bfx)}^\top\cdot \Sigma(\bfx)\bfxi(\bfx)\right)\geq \gamma\|\bfxi(\bfx)\|_{{\mathbb R}^3}^2. 
\end{equation*}
As before, this condition is satisfied if the eigenvalues of the matrix $\left(\begin{array}{rrcll} \tilde A & \tilde B \\ 
\overline{\tilde B} & \tilde C \end{array}\right)$ are  positive uniformly on $\Gamma$, where now 
$$\tilde A:=|\alpha|^2(\Re(\sigma_{11})\cos \theta+\Im(\sigma_{11})\sin \theta)), \qquad \tilde C:=|\beta|^2(\Re(\sigma_{22})\cos \theta+\Im(\sigma_{22})\sin \theta))$$
$$\tilde B:=\overline{\alpha}\beta\,\left(\frac{\sigma_{12}+\overline{\sigma}_{21}}{2}\cos \theta+\frac{\sigma_{12}-\overline{\sigma}_{21}}{2i}\sin \theta\right).$$
Thus the existence of the solution holds if for some $0\leq \theta\leq \pi/2$ and $\gamma>0$ constant  and for almost all $x\in \Gamma$ we have  
\begin{subequations}\label{existence}
\begin{eqnarray}
&\Re(\sigma_{11}+\sigma_{22})\cos \theta \geq \gamma,  \qquad \Im(\sigma_{11}+\sigma_{22})\sin \theta \geq \gamma, &\\[3pt]
& (\Re(\sigma_{11})\cos \theta+\Im(\sigma_{11})\sin \theta))(\Re(\sigma_{22})\cos \theta+\Im(\sigma_{22})\sin \theta))& \\
&\qquad \qquad \geq \displaystyle{\left|\frac{\sigma_{12}+\overline{\sigma}_{21}}{2}\cos \theta+\frac{\sigma_{12}-\overline{\sigma}_{21}}{2i}\sin \theta\right|^2}\nonumber.
\end{eqnarray}
\end{subequations}
Summarizing our requirements on $\Sigma$, throughout the paper we require that the surface tensor $\Sigma$ satisfies the following assumption which guarantees that the forward scattering problem  (\ref{eq:maxwell})-(\ref{eq:silver})  is well-posed, i.e. it has a unique solution in $X_{loc}(\curl, {\mathbb R}^3)$ depending continuously on the incident field.
\begin{assumption} \label{ass} The surface tensor $\Sigma\in L^\infty(\Gamma)^{2\times 2}$ satisfies conditions (\ref{unique}) and (\ref{existence}).
\end{assumption}
Note that Assumption \ref{ass} is quite general in that anisotropic surfaces are included in our analysis. If $\Re(\Sigma)$ is positive definite our assumptions include  the so-called highly directional hyperbolic meta-surfaces, for which the $\Im(\Sigma)$ is not sign-definite, i.e. has one positive and one negative eigenvalue at each point on $\Gamma$.  However, in the case of  resistive screens, i.e. when $\Re(\Sigma)\equiv 0$, we need $\Im(\Sigma)$ to be positive definite. Note also that we don't assume any symmetry on the tensor $\Sigma$ to possibly include  symmetry breaking  meta-surfaces (see e.g. \cite{bilow,Alu2, Alu,alu16,Alu3} and the references therein).

\section{The Inverse Scattering Problem} \label{inverse}

For an incident plane wave 
\[\bfE^{i}(\bfx; \bfd,\bfp):=\bfE^{i}(\bfx;\kappa, \bfd,\bfp)
\] given by (\ref{plane}) (since the wave number $\kappa$ is fixed  from now on we will drop the dependence of the fields on $\kappa$), the field far field pattern $\bfE_{\infty}(\hbfx;\bfd,\bfp)$  of the corresponding  scattered field is defined from the following asymptotic behavior of the scattered field {\cite{CK2019}}
\begin{align}\label{def_us}
\bfE^{s}(\hbfx;\bfd,\bfp)=\frac{\exp(i\kappa r)}{r}\left\{\bfE^{\infty}(\hbfx;\bfd,\bfp)+O\left(\frac{1}{r}\right)\right\}\mbox{ as } r:=|\bfx|\to\infty.
\end{align}
Our first goal is to prove a uniqueness theorem for the general inverse problem of determining $\Sigma$ from scattering data.  For this we need the following lemma, where  ${\mathbb{S}}:=\left\{\bfx\in {\mathbb R}^3: \, \Vert 
\bfx\Vert=1\right\}$  denotes the unit sphere in $\mathbb{R}^3$ :
\begin{lemma}\label{wind} Under Assumption \ref{ass}, the set 
$$\mbox{\em Span}\, \left\{\bfE_T(\cdot\,;\bfd,\bfp)|_{\Gamma}\quad \mbox{for all $\bfd\in {\mathbb S}$ and  $\bfp\in\mathbb{R}^3$, $\bfd\cdot\bfp=0$} \right\}$$
is dense in $\sL(\Gamma)$.
\end{lemma}
{\em Proof:} Assume that $\bfphi\in \sL(\Gamma)$ is such that
$$\int_{\Gamma}\bfphi \cdot {\bfE}_T(\cdot\,; \bfd,\bfp) \,dA={\mathbf 0}\qquad \mbox{for all $\bfd\in {\mathbb S}$ and  $\bfp\in\mathbb{R}^3$, $\bfd\cdot\bfp=0$}.$$
Let ${\mathbf U}\in X_{loc}(\curl, B_R)$ be the unique radiating solution (i.e. it satisfies  the Silver-M{\"u}ller radiation condition) of
\begin{eqnarray*}
\curl\curl \bfU-\kappa^2\bfU=0&&\qquad \mbox{ in }\mathbb{R}^3\setminus\Gamma\\
\bfU^+_T=\bfU^-_T&&\qquad \mbox{ on }\Gamma\\
\bfnu\times(\curl\bfU^+-\curl\bfU^-)-i\kappa\Sigma^{\top}\bfU^+_T =\bfphi && \qquad \mbox{ on }\Gamma.
\end{eqnarray*}
Note that  the transposed tensor $\Sigma^{T}$ satisfies  Assumption \ref{ass} since it does not involve any conjugation. 
Thus, noting that $\bfU^+=\bfU^-$ on $\Gamma$ and using the boundary condition for the total field $\bfE$,
\begin{eqnarray*}
 0&=&\int_{\Gamma}\left(\bfnu\times\curl\bfU^+-\bfnu\times\curl\bfU^--i\kappa\Sigma^\top \bfU_T\right)\cdot {\bfE}_T\, dA\\
&=& \int_{\Gamma}\left(\bfnu\times\curl\bfU^+-\bfnu\times\curl\bfU^-\right)\cdot \bfE_T-i\kappa \Sigma \bfE_T \cdot \bfU_T \, dA\\
&=&\int_{\Gamma}\left(\bfnu\times\curl\bfU^+-\bfnu\times\curl\bfU^-\right)\cdot \bfE_T-\left(\bfnu\times\curl\bfE^+-\bfnu\times\curl\bfE^-\right)\cdot \bfU_T\, dA\\
&=&\int_{\Gamma}\left(\bfnu\times\curl\bfU^+-\bfnu\times\curl\bfU^-\right)\cdot \bfE^s_T-\left(\bfnu\times\curl\bfE^{s+}-\bfnu\times\curl\bfE^{s-}\right)\cdot \bfU_T\, dA\\
&+&\int_{\Gamma}\left(\bfnu\times\curl\bfU^+-\bfnu\times\curl\bfU^-\right)\cdot \bfE^i_T-\left(\bfnu\times\curl\bfE^{i+}-\bfnu\times\curl\bfE^{i-}\right)\cdot \bfU_T\, dA.
\end{eqnarray*}
The first integral in the last sum is zero since both $\bfU$ and $\bfE^s$ are in $X_{loc}(\curl, B_R)$ (i.e their tangential traces across $\Gamma$ are continuous)  and are both radiating solutions to Maxwells equation. The second term in the second integral is also zero since $\curl \bfE^i$ doesn't jump  across $\Gamma$, but we keep it for use with integration by parts below. Thus noting that all jumps across $\partial D\setminus {\overline \Gamma}$ are zero, integrating by parts inside in $D$ and $B_R\setminus \overline D$, and using that $\bfU$ and $\bfE^i$ satisfy the same Maxwell's equations, we arrive at 
\begin{eqnarray*}
0&=&\int_{\Gamma}\left(\bfnu\times\curl\bfU^+-\bfnu\times\curl\bfU^-\right)\cdot \bfE^i_T-\left(\bfnu\times\curl\bfE^{i+}-\bfnu\times\curl\bfE^{i-}\right)\cdot \bfU_T\, dA\\
&=&\int_{\partial D}\bfnu\times\curl\bfU^+\cdot \bfE^i_T-\bfnu\times\curl\bfE^{i}\cdot \bfU_T\, dA\\
&-&\int_{\partial D}\bfnu\times\curl\bfU^-\cdot \bfE^i_T-\bfnu\times\curl\bfE^{i}\cdot \bfU_T\, dA\\
&=&\int_{B_R}\bfnu\times\curl\bfU\cdot \bfE^i_T-\bfnu\times\curl\bfE^{i}\cdot \bfU_T\, dA\\
&=&i\kappa \int_{\partial B_R}\left(\hat \bfx\times\curl\bfU(\bfx)\right)\cdot  (\bfd\times \bfp)\times \bfd e^{-i\kappa \bfd\cdot\bfx} +i\kappa \hat \bfx\times(\bfd\times \bfp)e^{-i\kappa \bfd\cdot\bfx} \cdot \bfU_T(\bfx)\, dA_{\bfx}
\end{eqnarray*}
for all $\bfd\in {\mathbb S}$ and  $\bfp\in\mathbb{R}^3$, $\bfd\cdot\bfp=0$, (note that $\bfp\exp(-i\kappa\bfd\cdot\bfx)$ is an incident field). Therefore  we have (see e.g.  \cite[Theorem 6.9]{CK2019})
$${\bf 0}=\bfd\times \int_{\partial B_R}\left[\frac{1}{i\kappa} \left(\hat \bfx\times\curl\bfU(\bfx)\right)\times \bfd+(\hat\bfx\times \bfU)\right]\cdot \bfp e^{-i\kappa \bfd\cdot\bfx}\,dA=\frac{4\pi}{i\kappa}\bfU^\infty(\hat{\bfx},\bfd)\cdot \bfp.$$
Since this holds  for all polarizations $\bfp$ we conclude that $\bfU^\infty=0$. Rellich's Lemma implies $\bfU={\bf 0}$ in ${\mathbb R}^3\setminus  \overline{\Gamma}$, whence $\bfphi={\bf 0}$ which concludes the proof.

\noindent
Now we are ready to prove a uniqueness theorem for the tensor $\Sigma$.
\begin{theorem}\label{U}
Assume that $\Sigma_1$ and $\Sigma_2$ satisfy Assumption \ref{ass} and that $\Gamma$ is a given piece-wise smooth open surface. Let  $\bfE^{\infty,1}(\hbfx;\bfd,\bfp)$ and $\bfE^{\infty,2}(\hbfx;\bfd,\bfp)$ be the far field pattern corresponding to the scattered fields $\bfE^{s,1}(\cdot;\bfd,\bfp)$ and $\bfE^{s,2}(\cdot;\bfd,\bfp)$ in $X_{loc}(\curl, {\mathbb R}^3)$ satisfying (\ref{eq:maxwell})-(\ref{eq:silver})  with $\Sigma_1$ and $\Sigma_2$ respectively, and incident plane wave $\bfE^{i}(\cdot;\bfd,\bfp)$ given by (\ref{plane}). If $\bfE^{\infty,1}(\cdot;\bfd,\bfp)=\bfE^{\infty,2}(\cdot;\bfd,\bfp)$ for all $\bfd\in {\mathbb S}$ and  $\bfp\in\mathbb{R}^3$ with $\bfd\cdot\bfp=0$, then $\Sigma_1=\Sigma_2$.
\end{theorem}
{\em Proof:}
Let $\bfU(\bfx):=\bfE^{s,1}(\hbfx;\bfd,\bfp)-\bfE^{s,2}(\hbfx;\bfd,\bfp)=\bfE^{1}(\hbfx;\bfd,\bfp)-\bfE^{2}(\hbfx;\bfd,\bfp)$. From the assumption we have $\bfU^\infty(\hat \bfx)={\bf 0}$ for $\hat \bfx\in {\mathbb S}$ and hence by Rellich Lemma $\bfU(\bfx)=0$ for all $\bfx\in {\mathbb R}^3\setminus\overline{\Gamma}$. Hence, noting that $\bfU_T={\bf 0}$, we have for almost all $\bfx \in \Gamma$ 
\begin{eqnarray*}
{\bf 0}&=&\bfnu\times(\curl\bfU^+-\curl\bfU^-)=i\kappa\Sigma_1\bfE_T^{1}(\hbfx;\bfd,\bfp)-i\kappa \Sigma_2\bfE_T^{2}(\hbfx;\bfd,\bfp)\\ 
&=& i\kappa(\Sigma_1-\Sigma_2)\bfE_T^{2}(\hbfx;\bfd,\bfp).
\end{eqnarray*}
Viewing $\Sigma_1-\Sigma_2$ as a linear operator on $L^2(\Gamma)$, the result follows from Lemma \ref{wind}.

\noindent
Note that the proof of Theorem \ref{U} shows that if $\Sigma$ is a piece-wise continuous scalar function, then  the far field pattern due to one incident plane waves uniquely determines it. Nevertheless, our target signatures   require the scattering data as stated in the next definition.

\begin{definition}[Inverse Problem]\label{scat}
{ The {inverse problem} we are concerned with is, provided that  the shape $\Gamma$ of the surface is 
known, determine indicators of changes in the surface tensor  $\Sigma$ from the scattering data. The {\em scattering data} is the set of 
the far field patterns $\bfE^{\infty}(\hbfx;\bfd,\bfp)\in L^2({\mathbb {S}})$  for all observation directions $\hbfx$ and incident directions 
$\bfd$ on the unit  sphere ${\mathbb{S}}$ and all $\bfp\in\mathbb{R}^3$, $\bfd\cdot\bfp=0$  at a fixed wave number $\kappa$.
}		
\end{definition}
{\em Remark:}
{\em It is important to emphasize that our theoretical study holds if the scattering data is given on a partial aperture, i.e. for observation directions $\hbfx\in {\mathbb S}_r\subset {\mathbb S}$ and incident directions 
$\bfd\in{\mathbb{S}}_t\subset {\mathbb S}$ and  two linearly independent polarization $\bfp$ such that $\bfp\cdot \bfd=0$, where receivers location ${\mathbb S}_r$ and transmitters locations ${\mathbb S}_t$ are open subsets (possibly the same) of the unit sphere. }

\bigskip

\noindent
The scattering data defines the {\it far field operator} $F: L_t^2({\mathbb{S}}) \rightarrow
L_t^2({\mathbb{S}})$  by
\begin{equation}\label{faropm}
(F\bfg)(\hat \bfx):=\int_{{\mathbb{S}}}\bfE^{\infty}(\hbfx;\bfd,\bfg(\bfd))ds_{\bfd},
\qquad \hat \bfx \in {\mathbb{S}}.
\end{equation}
Note that $F$ a linear operator since $\bfE^\infty$ depends linearly on polarization $\bfp$ by the linearity of the forward problem and linear dependence of the  incident wave on $\bfp$. It is bounded and  compact \cite{CBMS-kot}. By superposition $F\bfg$ is the
electric far field pattern of the scattered field solving   (\ref{eq:maxwell})-(\ref{eq:silver})  with $\bfE^i:=\bfE^i_\bfg$ where $\bfE^i_\bfg$ is the electric Herglotz wave function  with kernel $\bfg$ given by \cite[Section 6.6]{CK2019}
\begin{equation}\label{pairh}
\bfE^i_{\bfg}(\bfx)=i\kappa \int_{\bbS}e^{i\kappa \bfd\cdot \bfx}\bfg(\bfd)ds_{\bfd} \qquad g\in L_t^2(\bbS)
\end{equation}
which is an entire solution of the Maxwell's equations. A knowledge of the scattering data in  Definition \ref{scat}, implies a knowledge of the far field operator data. From now on the far field operator $F$ is the data for our target signatures. In the following we will denote by $\bfE_{\bfg}$, $\bfE^s_\bfg$ and $\bfE_\bfg^\infty$ the total electric field, the scattered electric field and the electric far field pattern, respectively, corresponding to   the electric Herglotz incident field $\bfE^i_\bfg$.

\bigskip

\noindent
Our target signatures are based  on a set of eigenvalues which can be determined from scattering data. This method makes use of a modification of the far field operator using an auxiliary impedance scattering problem, similar to that  introduced in \cite{Screens} for the Helmholtz equation. Given the particular features of Maxwell's system, we  adopt a slightly different approach to that used in~\cite{Screens} in order to avoid dealing with a mixed eigenvalue problem. Furthermore,  to restore the compactness  of the electromagnetic  Dirichlet-to-Neumann  operator,  we include a smoothing operator following \cite{CLM2017}. 

\noindent To this end we  recall the  linear operator ${\mathcal S}$  first introduced in   \cite{CLM2017, halla-steklov}:
\begin{equation}\label{Sdef0}
\begin{array}{rl}
\cS\,:\, H^{-1/2}(\curl_{\partial D},\partial D) & \longrightarrow \,H^{1/2}(\Div_{\partial D}^0,\partial D)\,\\
\bfv\ &\longmapsto \cS\bfv\,:=-\,\bfcurl_{\partial D} q\,, 
\end{array}
\end{equation}
where $q\in H^{1}(\partial D)/{\mathbb C}$ is the solution of the problem \[
\Delta_{\partial D} q=\curl_{\partial D} \bfv \mbox{ on }{\partial D}\] where $\Delta_{\partial D}$ is the surface Laplacian on $\partial D$
also given by $\Delta_{\partial D}q=\curl_{\partial D} \bfcurl_{\partial D}\,q$. In other words
for  $\bfv\in H^{-1/2}(\curl_{\partial D},\partial D)$ by
\begin{equation}
{\cal S}\bfv=-{\rm\mathbf{curl}}_{\partial D} \Delta_{\partial D}^{-1}{\rm{curl}}_{\partial D}\bfv\label{	}
\end{equation}
By using an eigensystem expansion (e.g. \cite{nedelec})  we see that $\bfcurl_{\partial D}\,q\in  H^{1/2}_t(\partial D)$. Thus, ${\mathcal S}\bfv\in H^{1/2}_t(\partial D)$, $\Div_{\partial D}\bfv=0$ and
\begin{equation*}
\|{\mathcal S}\bfv\|_{H^{1/2}(\Div_{\partial D}^0,\partial D)}= \|{\mathcal S}\bfv\|_{1/2,\partial D}=\|\curl_{\partial D} q\|_{1/2,\partial D}\leq C_{\mathcal S}\|\curl_{\partial D}\bfv\|_{-1/2,\partial D},
\end{equation*}
which means that  ${\mathcal S}$ is bounded linear operator. In addition, since $\curl_{\partial D} (\bfcurl_{\partial D} q-\bfv)=0$,  we can find $\varphi \in H^{1/2}(\partial B)$ such that $\bfcurl_{\partial D}q-\bfv=\nabla_{\partial D} \varphi$. Therefore, for all $\bfv\in H^{-1/2}(\curl_{\partial D},\partial D)$, there exist $q$ and $\varphi$ such that $\bfv\,=\,\bfcurl_{\partial D}q-\nabla_{\partial D} \varphi$, or, equivalently, ${\mathcal S}\bfv=\bfv+\nabla_{\partial D}  \varphi$.

\noindent 
We can now define the following auxiliary scattering problem for the field $\bfE^{(\lambda)}$:
\begin{subequations}\label{auxprob}
\begin{eqnarray}
\curl\curl \bfE^{(\lambda)}-\kappa^2\bfE^{(\lambda)}=0 && \quad\mbox{ in }\mathbb{R}^3\setminus\overline{D},\label{eq:maxwell_aux1}\\
\bfE^{(\lambda)}=\bfE^{(\lambda),s}+\bfE^{i}&& \quad \mbox{ in }\mathbb{R}^3\setminus D,\label{eq:scat_aux2}\\
\bfnu\times\curl\bfE^{(\lambda)}-\lambda {\cal S}\bfE_T^{(\lambda)}=0 && \quad \mbox{ on }\partial D,\label{eq:ban_aux3}\\\
\lim\limits_{|\bfx|\to \infty}\left(\curl \bfE^{{(\lambda)},s}\times \bfx-i\kappa|\bfx|\bfE^{{(\lambda)},s} \right)=0.&&\label{eq:silver_aux}
\end{eqnarray}
\end{subequations}
Here $\bfE^{(\lambda),s}$ denotes the scattered field for the above problem, and $\lambda\in\mathbb{C}$  is an auxiliary parameter which will play the role of the eigenvalue parameter used to find a target signature for $\Sigma$.  
\medskip

\noindent
To study the well-posedness of  (\ref{eq:maxwell_aux1})-(\ref{eq:silver_aux}) we recall from  \cite[Lemma 3.1]{CLM2017} that $\cS$ satisfies
\begin{equation}\label{propS:1}
\int_{\partial D} \cS\bfu_T\cdot\overline{\bfw_T}\,ds\,=\int_{\partial D} \bfu_T\cdot\overline{\cS\bfw_T}\,ds\,=\int_{\partial D} \cS\bfu_T\cdot\overline{\cS\bfw_T}\,ds\,,
\end{equation}
for all $\bfu,\,\bfw$ in $H(\curl, D)$ or $H(\curl, B_R\setminus {\overline{D}})$. Thus integrating by parts formally we have 
\begin{eqnarray}
&&\int_{B_R} (\curl \bfE^{(\lambda),s}\cdot\curl\overline{\bfv}-\kappa^2\bfE^{(\lambda),s}\cdot\overline{\bfv})\,dV-\lambda\int_{\partial D}\cS \bfE^s_T\cdot\overline{\bfv}_T\,dA \nonumber\\
&&\qquad \qquad +\int_{\partial B_R}\bfnu\times\curl\bfE^s\cdot \overline{\bfv}\,dA=\lambda\int_{\partial D}\cS\bfE^i_T\cdot\overline{\bfv}_T\,dA. \label{varEA}
\end{eqnarray}
From (\ref{propS:1}) by taking $\bfv:=\bfE^{(\lambda),s}$ and $\bfE^i={\bf 0}$ in (\ref{varEA}) in the same way as for the forward scattering problem we see that uniqueness is ensured if $\Im(\lambda)\geq 0$. Writing  $\int_{\partial B_R}\bfnu\times\curl\bfE^s\cdot \overline{\bfv}\,dA$ in terms of the exterior Calderon operator $G_e$ (c.f. \cite{monk}), we obtain  the existence of the solution  $\bfE^{(\lambda)}\in H_{loc}(\curl, \DE)$  by means of the Fredholm alternative \cite[Theorem 3.3]{CLM2017} stated in the theorem below.
\begin{theorem}
Assume that $\lambda\in\mathbb{C}$ is such that $\Im(\lambda)\geq 0$. Then the auxiliary problem (\ref{auxprob}) has a unique solution $\bfE^{(\lambda)}\in H_{loc}(\curl, \DE)$  depending continuously on the incident field $\bfE^i$.
\end{theorem}
\noindent
Let $\bfE^{(\lambda)}(\cdot;\bfd,\bfp)$ be the solution of (\ref{eq:maxwell_aux1})-(\ref{eq:silver_aux}) corresponding to the incident plane wave $\bfE^i:=\bfE^{i}(\cdot;\bfd,\bfp)$ and let $\bfE^{(\lambda),\infty}(\hbfx;\bfd,\bfp)\in L^2({\mathbb {S}})$ denote its far field pattern. The corresponding far field operator $F^{(\lambda)}: L_t^2({\mathbb{S}}) \rightarrow
L_t^2({\mathbb{S}})$ is
\begin{equation}\label{faroplambda}
(F^{(\lambda)}\bfg)(\hat \bfx):=\int_{{\mathbb{S}}}\bfE^{(\lambda),\infty}(\hbfx;\bfd, \bfg(\bfd))ds_{\bfd},
\qquad \hat \bfx \in {\mathbb{S}},
\end{equation}
which  is the far field pattern $\bfE_\bfg^{(\lambda),\infty}$ of the solution $\bfE_\bfg^{(\lambda), s}$ to (\ref{auxprob}) with incident field $\bfE^i:=\bfE^i_\bfg$  the electric Herglotz wave function  with kernel $\bfg$ given by (\ref{pairh}).

\medskip
\noindent
Next we define the {\it modified far field operator} ${\mathcal F}: L_t^2({\mathbb{S}}) \rightarrow L_t^2({\mathbb{S}})$ by
\begin{eqnarray}
({\mathcal F}\bfg)(\hat \bfx):&=&(F\bfg)(\hat \bfx)- (F^{(\lambda)}\bfg)(\hat \bfx) \label{modf}\\
&=&\int_{{\mathbb{S}}}\left[\bfE^{\infty}(\hbfx;\bfd, \bfg(\bfd))-\bfE^{(\lambda),\infty}(\hbfx;\bfd, \bfg(\bfd))\right]ds_{\bfd}.\nonumber
\end{eqnarray}
The study of injectivity of ${\mathcal F}$, allows us to arrive at an eigenvalue problem whose eigenvalues are the target signature for the thin screen. Indeed, assume  ${\mathcal F}\bfg={\bf 0}$, for some $\bfg\in \sL(\mathbb{S})$, $\bfg\not=0$, so that $\bfE_\bfg^{\infty}=\bfE_\bfg^{(\lambda),\infty}$ on $ {\mathbb S}$. By Rellich's lemma, $\bfE^s_\bfg=\bfE_\bfg^{(\lambda), s}$ in ${\mathbb R}^3\setminus\overline{D}$,  and the same holds true for the total fields $\bfE_\bfg=\bfE_\bfg^{(\lambda)}$. Using the boundary condition (\ref{eq:ban_aux3}) for $\bfE_\bfg^{(\lambda)}$ we obtain
$$\bfnu\times\curl\bfE_\bfg^{+}-\lambda {\cal S}\bfE_{\bfg T}^{^+}=0 \qquad \mbox{on}\; \partial D,$$
where again $+$ and $-$ indicate that we  approach the boundary from outside and inside, respectively.   On the other hand, from (\ref{trans1})-(\ref{trans2}) we have 
$$\bfE^+_{\bfg T}=\bfE^-_{\bfg T}\; \mbox{on } \partial D, \qquad  \bfnu\times\curl\bfE_{\bfg}^+=\bfnu\times \curl\bfE_{\bfg}^-\; \mbox{on } \partial D\setminus \Gamma, $$
$$\mbox{and} \qquad \bfnu\times\curl\bfE_\bfg^+=\bfnu\times \curl\bfE_\bfg^-+i\kappa\Sigma\bfE^+_{\bfg T}\; \mbox{on } \Gamma.$$
We can eliminate $\bfE^+_{\bfg T}$ using the above three relations, yielding the following  homogeneous problem for the total field $\bfE_g$ from inside $D$:
\begin{eqnarray*}
&\curl\curl \bfE_\bfg-\kappa^2\bfE_\bfg={\bf 0}&\qquad \mbox{ in } D,\\
&\bfnu\times\curl\bfE_\bfg+i\kappa\Sigma\bfE_T= \lambda {\cal S}\bfE_{\bfg T} & \qquad \mbox{ on }\Gamma,\\
&\bfnu\times\curl\bfE_\bfg= \lambda {\cal S}\bfE_{\bfg T} & \qquad \mbox{ on }  \partial D\setminus \Gamma.
\end{eqnarray*}
For fixed $\kappa$ we view this problem as an eigenvalue problem for $\lambda$. In particular, it is a modified Steklov type eigenvalue problem corresponding to the screen described by  ($\Gamma, \Sigma$). If this homogeneous problem has only the trivial solution, then $\bfE_\bfg={\bf 0}$ in $D$ and by continuity of the electromagnetic Cauchy data $\bfE_\bfg={\bf 0}$ in ${\mathbb R}^3\setminus {\Gamma}$. The jump conditions (\ref{trans1})-(\ref{trans2}) ensure that $\bfE_\bfg$ solves Maxwell's equations in ${\mathbb R}^3$  and, the fact that $\bfE_\bfg\equiv {\bf 0}$ implies that $\bfE^s_\bfg=-\bfE^i_\bfg$ in ${\mathbb R}^3$. Hence the Herglotz function  $\bfE^i_\bfg\equiv {\bf 0}$ as an entire solution of Maxwell's equations that satisfies the outgoing radiation condition,  whence $\bfg={\bf 0}$ (see e.g. \cite[Chapter 6]{CK2019}). 

\begin{definition}[$\Sigma$-Steklov Eigenvalues] \label{def-eig} Values of $\lambda\in {\mathbb C}$ with $\Im(\lambda)\geq 0$ for which 
\begin{subequations}\label{Sigma_eig}\begin{eqnarray}
&\curl\curl \bfw-\kappa^2\bfw={\bf 0}&\qquad \mbox{ in } D,\\
&\bfnu\times\curl\bfw+i\kappa\Sigma\bfw= \lambda {\cal S}\bfw_T & \qquad \mbox{ on }\Gamma,\\
&\bfnu\times\curl\bfw= \lambda {\cal S}\bfw_T & \qquad \mbox{ on }  \partial D\setminus \Gamma,
\end{eqnarray}
\end{subequations}
has non-trivial solution, are called $\Sigma$-Steklov eigenvalues.
\end{definition}
\noindent
We have proven the following result.
\begin{theorem}\label{inj}
Let $\Sigma$ satisfies Assumption \ref{ass}. If $\lambda$ is not  a $\Sigma$-Steklov eigenvalue, then the modified far field operator ${\mathcal F}: L_t^2({\mathbb{S}}) \rightarrow L_t^2({\mathbb{S}})$ is injective.
\end{theorem}
\noindent
Note that the converse is not true, i.e. if $\lambda$ is a $\Sigma$-Steklov eigenvalue this doesn't necessary imply that ${\mathcal F}$ is not injective.  Next we study the range of the compact modified far field operator. To this end we need to compute the $L^2$-adjoint  ${\mathcal F}_\Sigma^*$ adjoint of  the modified far field operator ${\mathcal F}_\Sigma$ corresponding  $\Sigma$.  
\begin{lemma}\label{reci}
 The adjoint ${\mathcal F}_{\Sigma}^*: L_t^2({\mathbb{S}}) \rightarrow L_t^2({\mathbb{S}})$ is given by
 $${\mathcal F}^*\bfg=\overline{R{\mathcal F}_{\Sigma^\top}R{\overline \bfg}}$$
where ${\mathcal F}_{\Sigma^\top}$ is the modified far field operator corresponding to the scattering problem (\ref{eq:maxwell})-(\ref{eq:silver}) with the coefficient  $\Sigma^\top$ (the transpose of the tensor $\Sigma$). Here $R: L_t^2({\mathbb{S}}) \rightarrow L_t^2({\mathbb{S}})$  is defined by $R\bfg(d):=g(-d)$.
\end{lemma}
{\em Proof:}
First,  in the same way as in the proof of  \cite[Theorem 6.30]{CK2019}, we can show that 
\begin{eqnarray*}
&& \qquad i\kappa 4 \pi \left\{\bfq\cdot \bfE^{(\lambda), \infty }(\hat\bfx; \bfd,\bfp)-\bfp\cdot \bfE^{(\lambda), \infty }(-\bfd; -\hat\bfx,\bfq) \right\} =\\
&&\int\limits_{\partial B_R}\left[\bf\nu\times  \bfE^{(\lambda)}(\cdot; \bfd,\bfp) \cdot \curl \bfE^{(\lambda)}(\cdot; -\hat\bfx,\bfq)-\bf\nu\times  \curl \bfE^{(\lambda)}(\cdot; \bfd,\bfp) \cdot \bfE^{(\lambda)}(\cdot; -\hat\bfx,\bfq)\right]\, dA\\
&& \qquad =0.
\end{eqnarray*}
\noindent
Then  using the boundary condition  (\ref{eq:ban_aux3}) and the fact that both fields satisfy the same Maxwell's equations in $B_R\setminus \overline{D}$ we obtain 
\begin{eqnarray}
&& \qquad i\kappa 4 \pi \left\{\bfq\cdot \bfE^{(\lambda), \infty }(\hat\bfx; \bfd,\bfp)-\bfp\cdot \bfE^{(\lambda), \infty }(-\bfd; -\hat\bfx,\bfq) \right\} \\\label{rec1}
&=\lambda &\int_{\partial D}\left[\bfE_T^{(\lambda)}(\cdot; \bfd,\bfp) \cdot {\mathcal S} \bfE_T^{(\lambda)}(\cdot; -\hat\bfx,\bfq)- {\mathcal S}\bfE_T^{(\lambda)}(\cdot; \bfd,\bfp) \cdot \bfE_T^{(\lambda)}(\cdot; -\hat\bfx,\bfq)\right]\, dA=0\nonumber
\end{eqnarray}
due to the symmetry of ${\mathcal S}$. Then,  the reciprocity relation 
$$\bfq\cdot \bfE^{(\lambda), \infty }(\hat\bfx; \bfd,\bfp)=\bfp\cdot \bfE^{(\lambda), \infty }(-\bfd; -\hat\bfx,\bfq), \; \mbox{ for all $\bfd$, $\hat \bfx$ in ${\mathbb S}$ and any two $\bfp, \bfq$ in ${\mathbb R}^3$}$$ 
used in the same way as in  \cite[Theorem 6.37]{CK2019}  shows that 
\begin{equation}\label{FL}
\left(F^{(\lambda)}\right)^*\bfg=\overline{R F^{(\lambda)}R{\overline \bfg}}.
\end{equation}
The above proof suggest that, since in general $\Sigma$ is not symmetric, to compute the adjoint $F^*_\Sigma$ we must consider the scattering problem with transpose  $\Sigma^\top$. Using arguments similar to the proof of (\ref{rec1}), we can prove
\begin{eqnarray*}
&& \qquad i\kappa 4 \pi \left\{\bfq\cdot {\bfE}_\Sigma^{(\lambda), \infty }(\hat\bfx; \bfd,\bfp)-\bfp\cdot {\bfE}_{\Sigma^\top}^{(\lambda), \infty }(-{\bfd}; -\hat{\bfx},{\bfq}) \right\} =\nonumber\\
&&\int\limits_{\partial B_R}\left[{\bf\nu}\times  {\bfE}_\Sigma^{(\lambda)}(\cdot; \bfd,\bfp) \cdot \curl {\bfE}_{\Sigma^\top}^{(\lambda)}(\cdot; -\hat\bfx,\bfq)-{\bf\nu}\times  \curl {\bfE}_\Sigma^{(\lambda)}(\cdot; \bfd,\bfp) \cdot {\bfE}_{\Sigma^\top}^{(\lambda)}(\cdot; -\hat\bfx,\bfq)\right]\, dA\nonumber\\
&&\qquad =0.\label{recf}
\end{eqnarray*}
where the subscript $\Sigma$ and $\Sigma^\top$ indicate that the fields correspond to the scattering problem (\ref{eq:maxwell})-(\ref{eq:silver}) with $\Sigma$ and $\Sigma^\top$, respectively. Again using the fact that both total fields solve the Maxwell's equation in $B_R\setminus{\Gamma}$ together with the jump conditions (\ref{trans1})-(\ref{trans2}) yield
\begin{eqnarray}
&& \qquad i\kappa 4 \pi \left\{\bfq\cdot \bfE_\Sigma^{(\lambda), \infty }(\hat\bfx; \bfd,\bfp)-\bfp\cdot \bfE_{\Sigma^\top}^{(\lambda), \infty }(-\bfd; -\hat\bfx,\bfq) \right\} \\\label{recf1}
&=&\int_{\Gamma}\left[\bfE_{\Sigma,T}^{(\lambda)}(\cdot; \bfd,\bfp) \cdot \Sigma^\top\bfE_{\Sigma^\top,T}^{(\lambda)}(\cdot; -\hat\bfx,\bfq)-  \Sigma{\bfE}_{\Sigma,T}^{(\lambda)}(\cdot; \bfd,\bfp) \cdot {\bfE}_{\Sigma^\top,T}^{(\lambda)}(\cdot; -\hat\bfx,\bfq)\right]\, dA=0.\nonumber
\end{eqnarray}
Then,  the reciprocity relation 
$$\bfq\cdot \bfE_\Sigma^{(\lambda), \infty }(\hat\bfx; \bfd,\bfp)=\bfp\cdot \bfE_{\Sigma^\top}^{(\lambda), \infty }(-\bfd; -\hat\bfx,\bfq), \;\mbox{ for all $\bfd$, $\hat \bfx$ in ${\mathbb S}$ and any two $\bfp, \bfq$ in ${\mathbb R}^3$}$$ 
now gives
\begin{equation}\label{FL2}
F_\Sigma^*\bfg=\overline{R F_{\Sigma^\top}R{\overline \bfg}}.
\end{equation}
Combining (\ref{FL}) and (\ref{FL2}) proves the result of the lemma.

\noindent
Lemma \ref{reci}  implies the following result about the range of  the modified far field operator ${\mathcal F}$. (Note that  in what follows  ${\mathcal F}$ denotes the modified operator corresponding to $\Sigma$.) 
\begin{theorem}\label{dense}
Let $\Sigma$ satisfies Assumption \ref{ass}. If $\lambda$ is not  a $\Sigma^\top$-Steklov eigenvalue, then the modified far field operator ${\mathcal F}: L_t^2({\mathbb{S}}) \rightarrow L_t^2({\mathbb{S}})$  has dense range.
\end{theorem}

\noindent
We close this section with some equivalent  expression related to  the operator ${\mathcal S}$, for later use.  From  \cite[Page 236]{CK2019} we have 
\[
{\curl}_{\partial D}\bfv=-\nabla_{\partial D}\cdot (\bfnu\times\bfv), 
\]
and since the vector surface curl  denoted ${\rm\mathbf{curl}}_{\partial D}$ is the adjoint of the scalar surface curl, we have
\[
{\rm\mathbf{curl}}_{\partial D}v=-\bfnu\times\nabla_{\partial D}v
\]
for a scalar function $v$ on $\partial D$.  We can then verify that
\[
{\rm curl}_{\partial D}{\rm\mathbf{curl}}_{\partial D}=-\Delta_{\partial D}.
\]
Using these relations we see that an equivalent definition of ${\cal S}$ is
\begin{equation}
{\cal S}\bfv=-\bfnu\times\nabla_{\partial D} \Delta_{\partial D}^{-1}\nabla_{\partial D}\cdot (\bfnu\times\bfv)\label{Sdefeq}
\end{equation}
and this is the expression we use  in our numerical experiments in Section \ref{numerics}. Note that
for any surface tangential vector $\bfv\in H^{-1/2}(\curl_{\partial D},\partial D)$
\[
{\rm curl}_{\partial D}({\cal S}\bfv-\bfv)=
(-{\rm curl}_{\partial D}{\rm\mathbf{curl}}_{\partial D} \Delta_{\partial D}^{-1}{\rm{curl}}_{\partial D}\bfv -{\rm curl}_{\partial D}\bfv)=0.
\]
From here we see that there exists a $v\in H^{1/2}(\partial D)$ such that 
\begin{equation}\label{decom}
{\cal S}\bfv=\bfv+\nabla_{\partial D}v.
\end{equation}
\section{The $\Sigma$-Steklov Eigenvalue Problem}\label{sigma-stek}
We can write the $\Sigma$-Steklov eigenvalue problem defined in Definition \ref{def-eig} in the equivalent variational form: Find $\bfw\in X(\curl, D)$ such that 
\begin{eqnarray}
&&\int_D\curl \bfw\cdot \curl\overline{\bfv}-\kappa^2\bfw\cdot \overline{\bfv}\, dV\label{evar}\\
&-&i\kappa \int_{\Gamma} \Sigma \bfw_T\cdot \overline{\bfv}_T\, dA+\lambda \int_{\partial D}{\cal S}\bfw_T\cdot {\cal S}\overline{\bfv}_T\,dA=0 \qquad \forall \,  \bfv\in X(\curl, D),\nonumber
\end{eqnarray}
where we have used (\ref{propS:1}) and recall that  the operator $\cS\,:\, H^{-1/2}(\curl_{\partial D},\partial D) \to \,H^{1/2}(\Div_{\partial D}^0,\partial D)$. 
\begin{prop} Let $\Sigma$ satisfy Assumption  \ref{ass}. 
\begin{enumerate}
\item If  $\Re\left(\overline{\bfxi(\bfx)}^\top \cdot \Sigma(\bfx) \bfxi(\bfx)\right)> 0$  a.e. $\bfx\in \Gamma$, $\forall  \bfxi$ tangential complex fields, then all $\Sigma$-Steklov eigenvalues $\lambda$ satisfy $\Im(\lambda)\geq 0$. Real eigenvalues $\lambda$ (if they exist) do not depend on $\Sigma$.
\item If $\Re(\Sigma)=0$ (the zero matrix) almost everywhere on $\Gamma$ then the eigenvalues maybe be real and complex. Complex eigenvalues appears in conjugate pairs.
\item If $\Re(\Sigma)=0$ (the zero matrix)  almost everywhere on $\Gamma$ and $\Im(\Sigma)$ is symmetric then the eigenvalue problem is self-adjoint hence all eigenvalues are real.
\end{enumerate}
\end{prop}
{\em Remark:} More generally  if $\Re\left(\overline{\bfxi}^\top \cdot \Sigma\bfxi\right)> 0$ in  $\Gamma_0\subseteq \Gamma$, the proof of Case {\em 1} shows that real eigenvalues (if they exists) do not carry information on $\Sigma$ in $\Gamma_0$

{\em Proof:} Suppose $\Im(\lambda)\leq 0$ and Case {\em 1} holds.  Letting $\bfv:=\bfw$ in (\ref{evar}) and taking the imaginary part, yields $\bfw_T=0$ on $\Gamma$. If $\Im(\lambda)< 0$ we obtain $\int_{\partial D}|{\cal S}\bfw_T|^2\,dA=0$  we obtain  $\cS\bfw_T={\bf 0}$ on $\partial D$ and from boundary condition also $\nu\times \curl\bfw={\bf 0 }$ on $\Gamma$. Hence $\bfw={\bf 0}$ in $D$ as a solution of the Maxwell's equation with zero Cauchy data on $\Gamma$. Furthermore,  real $\lambda$ are eigenvalues of the following  problem 
$$\curl\curl \bfw-\kappa^2\bfw={\bf 0}\quad \mbox{ in } D, \qquad \bfnu\times\curl\bfw= \lambda {\cal S}\bfw_T \quad \mbox{ on }  \partial D,$$ (which from  \cite{CLM2017} it has  an infinite sequence of real eigenvalues accumulating to $+\infty$) with corresponding eigenvectors satisfying $\bfw|_{\Gamma}=0$. Obviously, if they exists, do not depend on $\Sigma$. Case {\em 2} follows form the fact that all operators are real and it is sufficient to work on real Hilbert spaces. Case {\em 3} is obvious and is discussed later in this section.

\noindent
Using Helmholtz decomposition we have that
$$X(\curl, D)=X(\curl, \Div 0, D) \oplus \nabla P\qquad \mbox{where} \qquad P:=\left\{p\in H^1(D); \, p=0 \;\; \mbox{on}\; \partial D\right\} $$
$$\mbox{and} \qquad X(\curl, \Div 0, D):=\left\{\bfu\in X(\curl, D)\; \Div \bfu=0 \mbox{ in } D, \; \nu\cdot \bfu=0 \mbox {\; on } \partial D\setminus \Gamma \right\}.$$
We can split $\bfw=\bfw_0+\nabla w$, $\bfw_0\in X(\curl, \Div 0, D)$ and $ w\in P$. Using the fact that $\curl(\nabla w)=0$ and that $(\nabla w)_T=0$ and taking in (\ref{evar}) the test function  $\bfv=\nabla \xi$ for $\xi \in P$ we obtain  that $w$ satisfies $\int_D \nabla w \cdot \nabla \xi =0$, implying that $w=0$. Therefore we view (\ref{evar}) in $X(\curl, \Div 0, D)$.   By means of Riesz representation theorem, we define ${\mathbb A}_{\Sigma,\kappa}$, ${\mathbb T}_\kappa$, ${\mathbb S}:  X(\curl, \Div 0, D)\to X(\curl, \Div 0, D)$ by
\begin{equation*}\label{opA}
\left({\mathbb A}_{\Sigma,\kappa}\bfw, \bfv \right)_{X(\curl, D)}:=\int_D\curl \bfw\cdot \curl\overline{\bfv}+\bfw\cdot \overline{\bfv}\, dA-i\kappa \int_{\Gamma} \Sigma \bfw_T\cdot \overline{\bfv}_T\, dA,
\end{equation*}
\begin{equation*}\label{opT}
\left({\mathbb T}_{\kappa}\bfw, \bfv \right)_{X(\curl, D)}:=(\kappa^2-1)\int_D\bfw\cdot \overline{\bfv}\, dV,
\end{equation*}
\begin{equation*}\label{opS}
\left({\mathbb S}\bfw, \bfv \right)_{X(\curl, D)}:=\int_{\partial D}{\cal S}\bfw_T\cdot {\cal S}\overline{\bfv}_T\,dA=\int_{\partial D}{\cal S}\bfw_T\cdot \overline{\bfv}_T\,dA,
\end{equation*}
respectively. Then the  eigenvalue problem of finding the kernel of
$$({\mathbb A}_{\Sigma,\kappa}+{\mathbb T}_{\kappa}+\lambda{\mathbb S})\bfw={\bf 0}\qquad \qquad \bfw\in X(\curl, \Div0, D).$$
Since $\Sigma$ (not necessarily Hermitian) satisfies Assumption  \ref{ass} we have that the operator  (not necessarily selfadjoint) ${\mathbb A}_{\Sigma,\kappa}$ is coercive hence invertible. The selfadjoint operator ${\mathbb S}: X(\curl, \Div 0, D)\to X(\curl, \Div 0, D)$ is compact. Indeed let ${\bfw}_j\rightharpoonup  \bfw_0$ converges weakly to some  $\bfw_0\in X(\curl, \Div 0, D)$. By boundedness of the trace operator we have that $(\bfw_{j}-\bfw_{0})_T\rightharpoonup 0$ in $H^{-1/2}(\curl_{\partial D},\partial D)$ and by the boundedness of ${\mathcal S}$ we have ${\cal S}(\bfw_{j}-\bfw_{0})_T$ converges to $0$ weakly in $H^{1/2}(\Div_{\partial D}^0,\partial D)$ and strongly in $\sL(\partial D)$ by the compact embedding of the prior space to the latter.  Then
\begin{eqnarray*}
&&\|{\mathbb S}({\bfw}_j-{\bfw}_0)\|_{X(\curl, D)}^2=\int_{\partial D}{\cal S}(\bfw_{j}-\bfw_{0})_T\cdot   {\cS}\left(\overline{{\mathbb S}({\bfw}_j-{\bfw}_0)}\right)_T\,dA\\
&& = \int_{\partial D}{\cal S}(\bfw_{j}-\bfw_{0})_T\cdot  \left(\overline{{\mathbb S}({\bfw}_j-{\bfw}_0)}\right)_T\,dA\leq C\|{\cal S}(\bfw_{j}-\bfw_{0})_T\|_{\sL(\partial D)}\to 0 \;\mbox{strongly},
\end{eqnarray*}
where we use the trace theorem and the fact that $(\bfw_{j}-\bfw_{0})$ is bounded in $ X(\curl, \Div 0, D)$. The selfadjoint operator ${\mathbb T}_\kappa$ is also compact since  $X(\curl, \Div0, D)$ combined with the fact that $\nu \times \curl \bfu\in L^2(\partial D)$ and $\curl \bfu\in H(\curl, D)$,   is compactly embedded in $L^2(D)$  (see e.g. \cite{costabel}).  From the Analytic Fredholm Theory \cite{CK2019} we conclude that   ${\mathbb A}_{\Sigma,\kappa}+{\mathbb T}_{\kappa}+\lambda{\mathbb S}$ has non-trivial kernel for at most a discrete set of $\lambda\in {\mathbb C}$ without finite accumulation points, and is invertible with bounded inverse for $\lambda$ outside this set.

\noindent
From the above discussion, for the given wave number $\kappa$ we can choose a constant $\alpha$ such that for ${\bf f}\in H^{1/2}(\Div_{\partial D}^0,\partial D)$ the problem
\begin{subequations}\label{dtm}
\begin{eqnarray}
&\curl\curl \bfw-\kappa^2\bfw={\bf 0}&\qquad \mbox{ in } D,\label{dtn1a}\\
&\bfnu\times\curl\bfw+i\kappa\Sigma\bfw_T= \alpha {\cal S}\bfw_T +{\bf f} & \qquad \mbox{ on }\Gamma \label{dtn2a}\\
&\bfnu\times\curl\bfw=\alpha {\cal S}\bfw_T +{\bf f}& \qquad \mbox{ on }  \partial D\setminus \label{dtn3a}\Gamma.
\end{eqnarray}
\end{subequations}
has a unique solution in $X(\curl, D)$. Note that if $\Re(\overline{\bfxi}^\top \cdot \Sigma \bfxi)> 0$ on some open set $\Gamma_0\subseteq \Gamma$, one can choose $\alpha=0$. We define the operator ${\mathcal R}_\Sigma:H^{1/2}(\Div_{\partial D}^0,\partial D) \to H^{1/2}(\Div_{\partial D}^0,\partial D)$ mapping ${\bf f}\mapsto \cS\bfw_T$ where $\bfw$ solves (\ref{dtm}).
\begin{lemma}
${\mathcal R}_\Sigma:H^{1/2}(\Div_{\partial D}^0,\partial D) \to H^{1/2}(\Div_{\partial D}^0,\partial D)$ is a compact operator.
\end{lemma}
{\em Proof:} This Lemma is proven in \cite[Lemma 3.4]{CLM2017} for a slightly different problem. We include it  here for the reader convenience. Equation (\ref{dtn1a}) implies that $\curl \bfw\in H(\curl, \Div^0, D)$
and equations (\ref{dtn2a}) and (\ref{dtn3a}) imply that $\nu\times \curl \bfw \in \sL(\Gamma)$. From \cite{costabel} we conclude that $\bfw\in H^{1/2}(D)$ and $\nu\cdot \curl \bfw\in L^2(D)$ implying $\curl_{\partial D}\bfw_T=\nu\cdot \curl \bfw\in L^2({\partial D})$. But, by definition, there exists $q\in H^{1}(\partial D)/{\mathbb C}$ such that $\cS\bfw_T\,:=-\,\bfcurl_{\partial D} q\in H^{1/2}(\Div_{\partial D}^0,\partial D)$. Since $\curl_{\partial D}\bfcurl_{\partial D} q=\curl_{\partial D}\cS\bfw_T=\curl_{\partial D}\bfw_T\in L^2(\partial D)$ we obtain that $\bfcurl_{\partial D} q\in H^1_t(\partial D)$. Hence  $\cS\bfw_T\,:=-\,\bfcurl_{\partial D} q$ is in  $H^{1}(\Div_{\partial D}^0,\partial D)$.  The proof is completed by recalling the compact embedding of $H^{1}(\Div_{\partial D}^0,\partial D)$ into  $H^{1/2}(\Div_{\partial D}^0,\partial D)$.

\noindent
We have shown that $(\lambda, \bfw)$ is an eigen-pair of the $\Sigma$-Steklov eigenvalue problem if and only if $\left(\frac{1}{\lambda-\alpha}, \cS\bfw_T\right)$ is an eigenpair of  the compact operator ${\mathcal R}_\Sigma$.
\begin{lemma} \label{l4}
 Let $\Sigma^\top$ be the transpose of $\Sigma$. If $\lambda$ is  a $\Sigma^\top$-Steklov eigenvalue then $1/(\lambda-\alpha)$ is an eigenvalue of  ${\mathcal R}_{\Sigma^\top}:H^{1/2}(\Div_{\partial D}^0,\partial D) \to H^{1/2}(\Div_{\partial D}^0,\partial D)$ which maps $h\mapsto {\cS}\bfv_T$ where $\bfv\in X(\curl, D)$ solves
\begin{subequations}\label{dtm2}
\begin{eqnarray}
&\curl\curl \bfv-\kappa^2\bfv={\bf 0}&\qquad \mbox{ in } D,\label{dtn1}\\
&\bfnu\times\curl\bfv+i\kappa\Sigma^\top\bfv_T= \alpha {\cal S}\bfv_T +{\bf h}& \qquad \mbox{ on }\Gamma \label{dtn2}\\
&\bfnu\times\curl\bfv=\alpha {\cal S}\bfv_T +{\bf h}& \qquad \mbox{ on }  \partial D\setminus \label{dtn3}\Gamma.
\end{eqnarray}
\end{subequations}
Furthermore ${\mathcal R}_{\Sigma^\top}$ is the transpose  (Banach adjoint) operator ${\mathcal R}^{\top}_\Sigma$  of ${\mathcal R}_\Sigma$,  where we have identified the Sobolev space $H^{1/2}(\Div_{\partial D}^0,\partial D)$ with its dual. In particular the set of $\Sigma^\top$-Steklov eigenvalues coincides with the set of $\Sigma$-Steklov eigenvalues.

\end{lemma}

{\em Proof:}
First note that if $\Sigma$ satisfies Assumption  \ref{ass} so does $\Sigma^\top$, hence the characterization of $\Sigma^\top$-Steklov eigenvalues follows form the above discussion. Next, let ${\bf f},{\bf h}\in H^{1/2}(\Div_{\partial D}^0,\partial D)$ and $\bfw$ and $\bfv$ such that ${\mathcal R}_\Sigma{\bf f}=\cS\bfw_T$ and ${\mathcal R}_{\Sigma^\top} {\bf h}=\cS\bfv_T$, where $\bfw$ and $\bfv$ satisfy (\ref{dtm}) and (\ref{dtm2}), respectively. Then we have
\begin{eqnarray*}
0&=&\int_D\curl \bfw\cdot \curl{\bfv}-\kappa^2\bfw\cdot {\bfv}\, dV\label{evara}\\
&-&i\kappa \int_{\Gamma} \Sigma \bfw_T\cdot {\bfv}_T\, dA+\alpha\int_{\partial D}{\cal S}\bfw_T\cdot {\cal S}{\bfv}_T\,dA +\int_{\partial D}{\bf f}\cdot {\cal S}{\bfv}_T\,dA
\end{eqnarray*}
and
\begin{eqnarray*}
0&=&\int_D\curl \bfv\cdot \curl{\bfw}-\kappa^2\bfv\cdot {\bfw}\, dV\label{evarb}\\
&-&i\kappa \int_{\Gamma} \Sigma^\top \bfv_T\cdot {\bfw}_T\, dA+\alpha\int_{\partial D}{\cal S}\bfv_T\cdot {\cal S}{\bfw}_T\,dA +\int_{\partial D}{\bf h}\cdot {\cal S}{\bfw}_T\,dA.
\end{eqnarray*}
where we have used (\ref{decom}), the fact that $\Div_{\partial D}{\bf f}=\Div_{\partial D}{\bf h}=0$ and the Helmholtz orthogonal decomposition $\bfmu=\mathbf{curl}_{\partial D} q+\nabla_{\partial D}p$ for any tangential field $\bfmu$ on the boundary. The above yields
$$\int_{\partial D}{\bf f}\cdot {\cal S}{\bfv}_T\,dA=\int_{\partial D}{\bf h}\cdot {\cal S}{\bfw}_T\,dA.$$
This proves that ${\mathcal R}_\Sigma^\top={\mathcal R}_{\Sigma^\top}$. The fact that they have the same non-zero eigenvalues follows for the Fredholm theory for compact operators, more precisely that  for $\eta\neq 0$, the dimension of $\mbox{Kern}({\mathcal R}_{\Sigma} -\eta I)$ and  $\mbox{Kern}({\mathcal R}_{\Sigma}^\top -\eta I)$ coincide.

\noindent
Thus we have shown that if $\Sigma$ satisfies  Assumption  \ref{ass} then the set of $\Sigma$-Steklov eigenvalues is discrete without finite accumulation points. The existence of (possibly complex) $\Sigma$-Steklov eigenvalues  could be proven by adapting the approach in \cite{halla-steklov}. We don't pursue this investigation here since it is out of the scope of the paper. 

\medskip

{\bf{The self-adjoint case}}. If  $\Sigma$ is  symmetric and $\Re(\Sigma)=0$ a.e. in $\Gamma$,   then ${\mathcal R}_{\Sigma}$ is compact and self-adjoint. Note that Assumption  \ref{ass} implies that $\Im(\Sigma)$ is positive definite.  In this case  $\Sigma$-Steklov eigenvalues  $\{\lambda_j\}$ form an infinite sequence of real numbers without finite accumulation point. We have seen that $\mu_j=\frac{1}{\lambda_j-\alpha}$, where  $\{\mu_j, \bfphi_j\}$ is an eigenpair of the compact self-adjoint operator ${\mathcal R}_{\Sigma}$, and  that by Hilbert-Schmidt theorem the eigenfunctions $\bfphi_j$ form a orthonormal basis for  $H^{1/2}(\Div_{\partial D}^0,\partial D)$.  To obtain additional estimates in this case we need the assumption
\begin{assumption}\label{ass2}
The wave number $\kappa$ is such that the homogeneous problem
\begin{eqnarray*}
 &\curl \bfw\curl \bfw-\kappa^2 \bfw={\bf{0}} \quad  \mbox{in}\, D&\\
&\bfnu \times \curl \bfw={\bf0} \quad  \mbox{on \; $\partial D\setminus \overline{\Gamma}$}\qquad \bfnu \times \curl \bfw=\Im(\Sigma) \bfw_T \quad \mbox{on\; $ \Gamma$}&
\end{eqnarray*}
has only the trivial solution.
\end{assumption}   
\begin{theorem}
Under Assumption \ref{ass2} there are finitely many positive  $\Sigma$-Steklov eigenvalues, thus the eigenvalues accumulate to $-\infty$.

\end{theorem}
{\em Proof:}Assume to the contrary  that there exists a sequence of distinct $\lambda_j>0$ converging to $\infty$. Denote by $\bfw_j$ the  solution of (\ref{dtm}) in $X(\curl, D)$ with ${\bf f}:=\bfphi_j$. We may normalize the sequence  $\|\bfw_j\|_{X(\curl, D)}+\|\bfw_{j,T}\|_{L^2(\partial D)}=1$.
 Furthermore since $(\lambda_j-\alpha){\mathcal S}\bfw_{j,T}=(\lambda_j-\alpha){\mathcal R}_{\Sigma}\bfphi_j=\bfphi_j$ we have
\begin{eqnarray*}
&&\int_D|\curl \bfw_j|^2-\kappa^2|\bfw_j|^2dV+\kappa \int_{\Gamma} \Im(\Sigma) \bfw_{j,T}\cdot {\bfw}_{j,T}\, dA+\alpha\int_{\partial D}{\cal S}\bfw_{j,T}\cdot {\bfw}_{j,T}\,dA\nonumber\\
&&\qquad \qquad  =(\alpha-\lambda_j)\int_{\partial D} {\cal S}{\bfw}_{j,T}\cdot \bfw_{j,T}\,dA
\end{eqnarray*}
which from  (\ref{propS:1})  gives
\begin{equation}\label{varas}
\int_D|\curl \bfw_j|^2-\kappa^2|\bfw_j|^2dV+\kappa \int_{\Gamma} \Im(\Sigma) \bfw_{j,T}\cdot {\bfw}_{j,T}\, dA=-\lambda_j\int_{\partial D} |{\cal S}{\bfw}_{j,T}|^2\,dA.
\end{equation}
Since the left-hand side is bounded we conclude that ${\cal S}{\bfw}_{j,T}\to 0$ in $L^2(\partial D)$ as $j\to \infty$. Next, a subsequence of $\bfw_j$ converges weakly to some $\bfw\in X(\curl, D)$. Since for all ${\bf z}\in X(\curl, D)$ we have 
$$\int_D\curl \bfw_j\cdot \curl {\bf z}-\kappa^2\bfw_j\cdot{\bf z}\,dV+\kappa \int_{\Gamma} \Im(\Sigma) \bfw_{j,T}\cdot {\bf z}_T\, dA=-\lambda_j\int_{\partial D} {\cal S}{\bfw}_{j,T} \cdot {\bf z}_T\,dA$$
we conclude that  the weak limit  satisfies the problem in Assumption \ref{ass2}, thus $\bfw ={\bf0}$. Using the Helmholtz decomposition and noting that $\Div \bfw_j=0$ and  $\kappa^2 \bfnu\cdot\bfw_j= \nu\times \curl \bfw_j\in L^2(\partial D)$  we conclude that $\bfw_j \rightharpoonup {\bf0}$ in $H^{1/2}(D)$ hence $\bfw_j\to {\bf0}$ strongly in $L^2(D)$. From (\ref{varas}) since $\Im(\Sigma)$ is positive and all $\lambda_j>0$ we have that 
$$\int_D|\curl \bfw_j|^2-\kappa^2|\bfw_j|^2dV+\kappa \int_{\Gamma} \Im(\Sigma) \bfw_{j,T}\cdot {\bfw}_{j,T}\, dA<0,$$
thus $\curl \bfw_j\to \bf0$ is $L^2(D)$ and  $\bfw_{j,T}\to \bf0$ in $L^2(\Gamma)$ contradicting the normalization.

\medskip
\noindent
The above discussion suggests that  if Assumption \ref{ass2} is satisfied, $\alpha>0$ can be chosen large enough such that all eigenvalues of ${\mathcal R}_{\Sigma}$ are negative. Using the Fischer-Courant max-min principle applied to the positive compact self-adjoint operator $-{\mathcal R}_\Sigma$, we have 
$$\mu_j=\max\limits_{U_{j-1}\in{\mathcal U}_{j-1}}\min\limits_{{\bf {f}}\in U_j, {\bf f}\neq {\bf 0}} \frac{\left({\mathcal R}_\Sigma{\bf f}, {\bf f}\right)_{H^{1/2}(\Div_{\partial D}^0,\partial D)}}{\|{\bf f}\|_{H^{1/2}(\Div_{\partial D}^0,\partial D)}^2}$$
where ${\mathcal U}_\ell$ is the set of all  linear subspace of $H^{1/2}(\Div_{\partial D}^0,\partial D)$ of  dimension $\ell$, $\ell= 1,2 \cdots$, which can be used to understand monotonicity of $\Sigma$-Steklov eigenvalues in terms of surface tensor $\Sigma$.


\section{Numerical Solution of the Inverse Problem}\label{numerics}
We propose a solution method for  the inverse problem formulated in Definition \ref{scat}. This method is based on a target signature that is computable from the scattering data defined  in Definition \ref{scat}. The target signature is defined precisely below.

\begin{definition}\label{TS}{\em [Target Signature for the Surface Tensor $\Sigma$]} Given $\Gamma$ piece-wise smooth and a domain $D$ with $\Gamma\subset \partial D$ the target signature for the unknown  surface tensor $\Sigma$ that satisfies Assumption \ref{ass}, is the set of $\Sigma$-Steklov eigenvalues defined in Definition \ref{def-eig}.
\end{definition}
This section is devoted to a discussion on  how the target signature is determined from the scattering and presenting numerical experiments showing the viability of our approach.  But, before providing preliminary numerical examples to illustrate our theory, we first give some general details about the results. Four pieces of software are needed for this purpose which we describe next.  All finite element implementations were performed using NGSolve~\cite{netgen}.
\subsection{Synthetic scattering data}
We need to find ${\cal F}$ which in turn requires solving the forward and auxiliary-forward problem as follows:
\begin{enumerate}
\item We use synthetic (computed) far field data so we need to approximate the forward problem (\ref{Forward_problem}). This is accomplished either using a standard edge finite element solver with a Perfectly Matched Layer (PML) to terminate the computational region.
\item We need to solve the auxiliary forward problem (\ref{auxprob}) for many choices of the parameter $\lambda$. This is done using edge finite elements and the PML.
\end{enumerate}

\subsection{Determination of $\Sigma$-Steklov eigenvalues from scattering data}\label{det}
We start by  discussing  the theoretical framework for the determination of $\Sigma$-Steklov eigenvalues from a knowledge of the modified far field operator ${\mathcal F}$. Note that ${\mathcal F}=F-F^{(\lambda)}$ is available to us since $F$  is known from the measured scattering data, whereas $F^{(\lambda)}$ for given $\Gamma$,  is computed by solving the auxiliary problem (\ref{auxprob}) which does not involve the unknown $\Sigma$. Note that, in practice, when problems of nondestructive testing of thin inhomogeneities, $F^{(\lambda)}$ can be precomputed and stored for a set of $\lambda \in {\mathbb C}$, $\Im(\lambda)\leq 0$, and this set may possibly  be determined using  a-priori information on the electromagnetic material properties encoded in $\Sigma$.   

\noindent
In view of  Theorem \ref{dense} and Lemma \ref{l4} we  now have the following result which is the fundamental theoretical  ingredient if the determination of $\Sigma$-eigenvalues from scattering data.
\begin{theorem}\label{dense2}
Let $\Sigma$ satisfy Assumption \ref{ass}. If $\lambda\in {\mathbb C}$ is not  a $\Sigma$-Steklov eigenvalue, then the modified far field operator ${\mathcal F}: L_t^2({\mathbb{S}}) \rightarrow L_t^2({\mathbb{S}})$  is injective and  has dense range.
\end{theorem}
\noindent
Using Theorem \ref{dense2}, an appropriate factorization ${\mathcal F}$ along with a denseness property of the total fields $\bfE_\bfg^{(\lambda)}$ solutions  to (\ref{auxprob}) with incident field $\bfE^i:=\bfE^i_\bfg$  the Herglotz wave function and finally making use of the Fredholm property of the resolvent of the $\Sigma$-Steklov eigenvalue problem it is possible to show the following result. To avoid repetition, for the proof of this result,  we refer the reader to \cite{heejin22} for the same problem but in the scalar case, to \cite{CLM2017} for a slightly different  problem but for the vectorial Maxwell's equations, and to \cite{CCH16} for a comprehensive discussion of this matter.  
Let ${\bf E}_{e,\infty}(\hat{\bf x},{\bf z},{\bfq})$ denote the far field pattern of the electric dipole with source at ${\bf z}$ and with polarization ${\bf q}$ given by
\[
{\bf E}_{e,\infty}(\hat{\bf x},{\bf z},{\bfq})=\frac{i\kappa}{4\pi} (\hat{\bf x}\times{\bf q})\times \hat{\bf x}\exp(-i\kappa \hat{\bf x}\cdot{\bf z}).
\]

\begin{theorem}\label{determination}
Let $\Sigma$ satisfy Assumption \ref{ass} and $\Gamma$ be a piece-wise smooth open surface embedded in a closed surface $\partial D$ circumscribing a connected region $D$. The following dichotomy  holds: 
\begin{itemize}
\item[(i)] Assume that $\lambda\in {\mathbb C}$ is not  a $\Sigma$-Steklov eigenvalue, and $z\in D$. Then there exists a sequence $\left\{{\bf g}^z_n\right\}_{n\in {\mathbb N}}$ in $L_t^2({\mathbb{S}}) $ such that
\begin{equation}\label{57}
\lim_{n\to 0}\|{\cal F}{\bf g}^z_n(\hat{\bf x})- {\bf E}_{e,\infty}(\hat{\bf x},{\bf z},{\bfq})\|_{L_t^2({\mathbb{S}})}=0
\end{equation}
and $\|\bfE_{{\bf g}^z_n}\|_{X(\curl, D)}$ remains bounded.
\item[(ii)]  (i) Assume that $\lambda\in {\mathbb C}$ is a $\Sigma$-Steklov eigenvalue. Then, for every  sequence $\left\{{\bf g}^z_n\right\}_{n\in {\mathbb N}}$
satisfying (\ref{57}), $\|\bfE_{{\bf g}^z_n}\|_{X(\curl, D)}$ cannot be bounded for any $z\in D$, except for a nowhere dense
set.
\end{itemize}
\end{theorem}
This theorem suggest that an ``approximate" solution ${\bf g}\in \sL(\mathbb{S}^2)$ of  the first kind integral equation 
\begin{equation}
{\cal F}{\bf g}(\hat{\bf x})= {\bf E}_{e,\infty}(\hat{\bf x},{\bf z},{\bfq})\mbox{ for all }\hat{\bf x}\in \mathbb{S}, \mbox{ and $z\in D$}
\label{far-field-eq}
\end{equation}
becomes unbounded if $\lambda\in {\mathbb C}$ hits  a $\Sigma$-Steklov eigenvalue. We remark that the procedure of  computing   $\left\{{\bf g}^z_n\right\}_{n\in {\mathbb N}}$ with the particular  behavior explained in Theorem \ref{determination}, can be made rigorous by applying the so-called generalized linear sampling method \cite[Chapter 5]{CCH16}.   Equation (\ref{far-field-eq})  is ill-posed since ${\cal F}$ is compact, but can be solved approximately using Tikhonov regularization for any choice of ${\bf z}$ and ${\bf q}$.   For the calculation of target signatures, we discretize (\ref{far-field-eq}) using the incident directions as quadrature points on $\partial D$, and chose $\hat{\bf x}$ to be the measurement points. In the results to be presented here we use 96 incoming plane wave directions and the same number of measurement points and assume that the polarization and phase of the far field pattern is available at each measurement point.  Then assuming that $D$ is a priori known, we take several random choices of ${\bf z}\in D$ (15 in our examples below).  For each point, and for the three canonical polarizations we solve the far field equation (\ref{far-field-eq}) approximately using Tikhonov regularization and average the norms of the three resulting ${\bf g}$ for the random points ${\bf z}$.  This is solved for a discrete choice of $\lambda$ in the 
interval in which it is desired to detect eigenvalues.  Peaks in the averaged norm of ${\bf g}$ are expected to coincide with $\Sigma$-Steklov eigenvalues.

\subsection{Direct calculation of $\Sigma$-Steklov eigenvalues} \label{eigen-num} To check the performance of our method for identifying $\Sigma$-Steklov eigenvalues, we also need to approximate the eigenvalue problem (\ref{Sigma_eig}) and this is again accomplished using finite elements.
  For $\bfw\in X({\rm curl},D)$ we introduce an auxiliary variable $z\in H^1(\partial D)/\mathbb{C}$ that satisfies
\[
\Delta_{\partial D}z=\nabla_{\partial D}\cdot(\bfnu\times \bfw)
\]
so ${\cal S}\bfw=-\bfnu\times\nabla_{\partial D}z$. We rewrite \eqref{Sigma_eig} as the problem of finding $z\in H^1(D)/\mathbb{C}$ and non-trivial
$\bfw\in H(\curl;D)$ and $\lambda\in \mathbb{C}$  such that
\begin{subequations}\label{eigprob1}
\begin{eqnarray}
\curl\curl \bfw-\kappa^2\bfw&=&0\mbox{ in }D,\label{eq:maxwell_eig1}\\
\bfnu\times\curl\bfw+i\kappa\Sigma\bfw_T&=&-\lambda \bfnu\times\nabla_{\partial D}z\mbox{ on }\Gamma,\\
\bfnu\times\curl\bfw&=&-\lambda \bfnu\times \nabla_{\partial D} z\mbox{ on }\partial D\setminus\Gamma,\\
\Delta_{\partial D} z-\nabla_{\partial D}\cdot(\bfnu\times \bfw)&=&0\mbox{ on }\partial D.
\end{eqnarray}
\end{subequations}
 Multiplying \eqref{eq:maxwell_eig1} by the complex conjugate of a test function $\bfv\in X({\rm curl};D)$, integrating by parts and using the boundary conditions in (\ref{eigprob1}), we obtain:
\begin{eqnarray*}
&&\int_{D} (\curl \bfw\cdot\curl\overline{\bfv}-\kappa^2\bfw\cdot\overline{\bfv})\,dV -\lambda\int_{\partial D}\bfnu\times\nabla_{\partial D}z\cdot\overline{\bfv}_T\,dA\\&&\quad-i\kappa\Sigma\int_{  \Gamma}\bfw_T\cdot\overline{\bfv}_T\,dA=0.
\end{eqnarray*}
So we define $A^{\rm eig}, b^{\rm eig}:(X({\rm curl},D)\times H^1(D)\times \mathbb{C})\times (X({\rm curl},D)\times H^1(D)\times \mathbb{C}) \to \mathbb{C}$ by
\begin{eqnarray*}
&&a^{\rm eig}((\bfw,z,r),(\bfv,q,s))
=\int_{D} (\curl \bfw\cdot\curl\overline{\bfv}-\kappa^2\bfw\cdot\overline{\bfv})\,dV -i\kappa\Sigma\int_{  \Gamma}\bfw_T\cdot\overline{\bfv}_T\,dA\\
&&\quad \qquad  +\int_{\partial D}\nabla_{\partial D} z\cdot\nabla_{\partial D}\overline{q}\, dA -{\int_{\partial D}\bfnu\times\bfw\cdot\nabla_{\partial D}\overline{q}\,dA}+\int_{\partial D}z\overline{s}-\overline{q}r\,dA\\
&&b^{\rm eig}((\bfw,z,r),(\bfv,q,s))={\int_{\partial D}\bfnu\times\nabla_{\partial D}z\cdot\overline{\bfv}_T\,dA}
\end{eqnarray*}
and seek non-trivial $(\bfw,z,r)\in X({\rm curl},D)\times H^1(D)\times \mathbb{C}$ and $\lambda\in \mathbb{C}$ such that
\[
a^{\rm eig}((\bfw,z,r),(\bfv,q,s))=\lambda b^{\rm eig}((\bfw,z,r),(\bfv,q,s)),
\]
for all $(\bfv,q,s)\in X({\rm curl},D)\times H^1(D)\times \mathbb{C}$. This can be discretized using edge and vertex finite elements. 
\subsection{Examples} 
\begin{figure}
\begin{center}
\resizebox{0.45\textwidth}{!}{\includegraphics{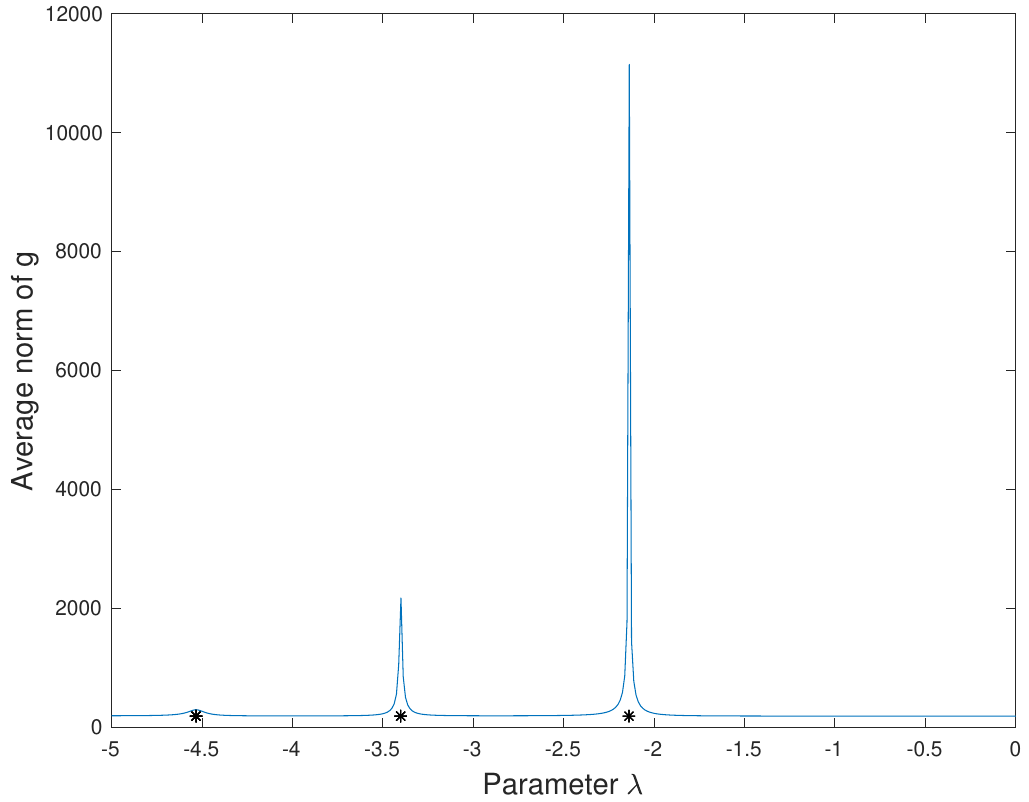}}
\end{center}
\caption{Target signatures for the full unit sphere at $\kappa=1.9$ and $\Sigma=(0.5i)I$.  We show results computed from the far field pattern as the curve of the average norm of ${\bf g}$ against the auxiliary parameter $\eta$.  We also show the first three $\Sigma$-Steklov eigenvaues  marked as $*$. Peaks of the avergae norm of ${\bf g}$ correspond well to $\Sigma$-Steklov eigenvalues.}\label{sphere}
\end{figure}

\begin{figure}
\begin{center}
\resizebox{0.45\textwidth}{!}{\includegraphics{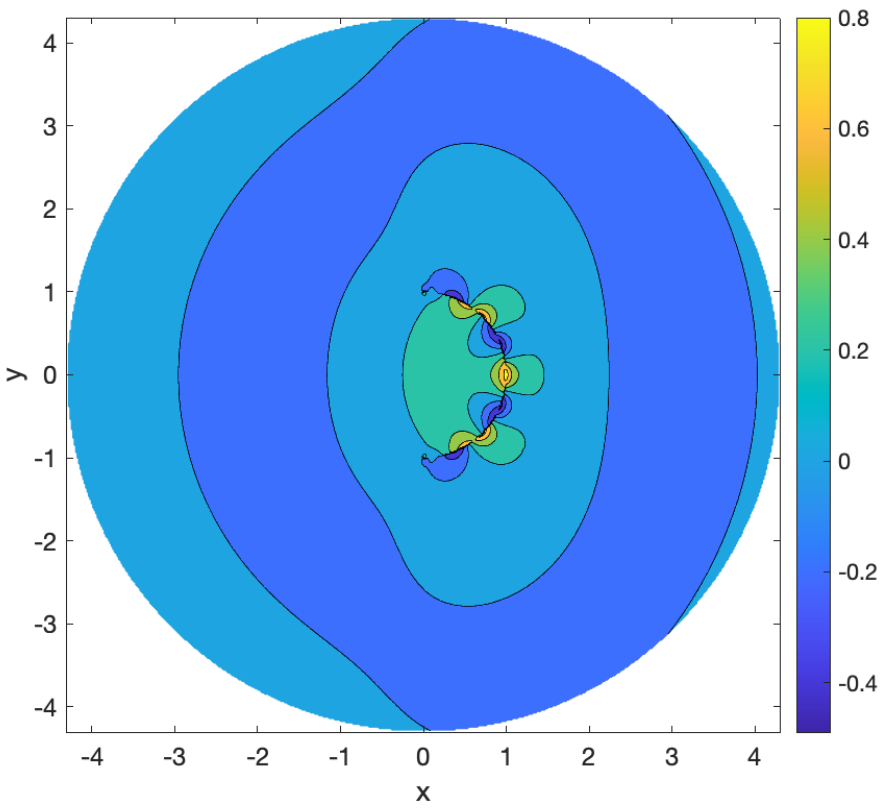}}
\end{center}
\label{fwd}
\caption{A contour map of the real part of the third component of the scattered electric field in the plane $z=0$.  Creeping waves along the screen are clearly visible.  These waves have a shorter wavelength than the field in the bulk,
so imposing an additional computational burden on the forward solver.}
\end{figure}

\begin{figure}
\begin{center}
\resizebox{0.45\textwidth}{!}{\includegraphics{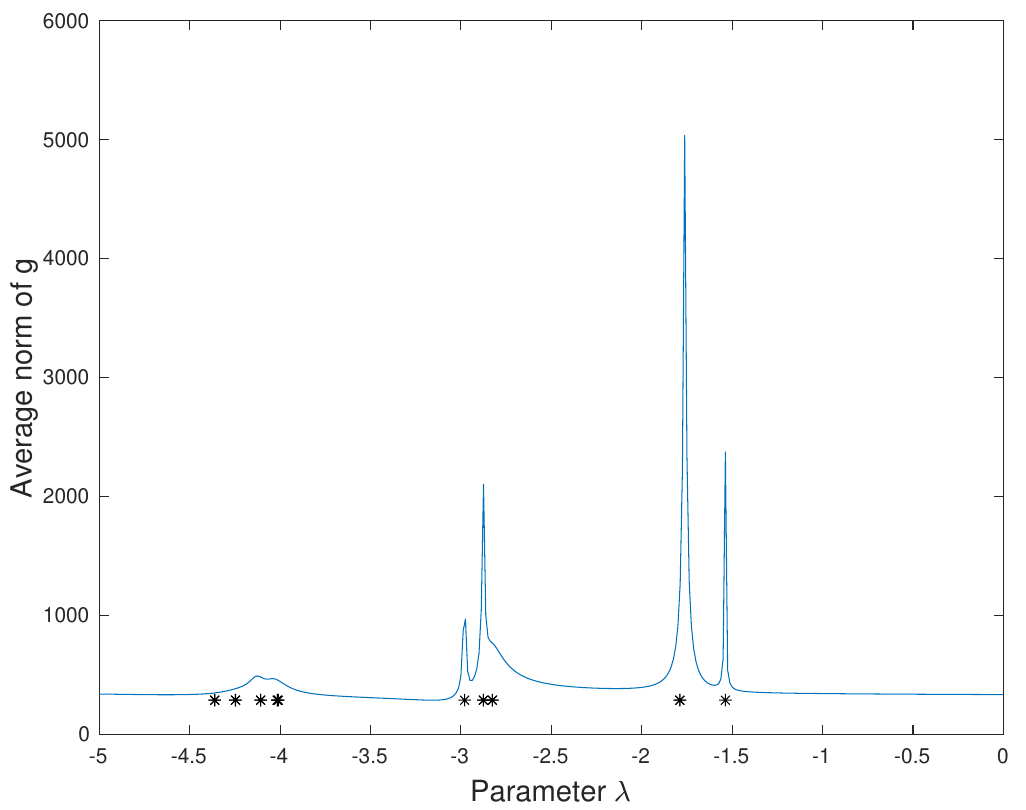}}
\resizebox{0.45\textwidth}{!}{\includegraphics{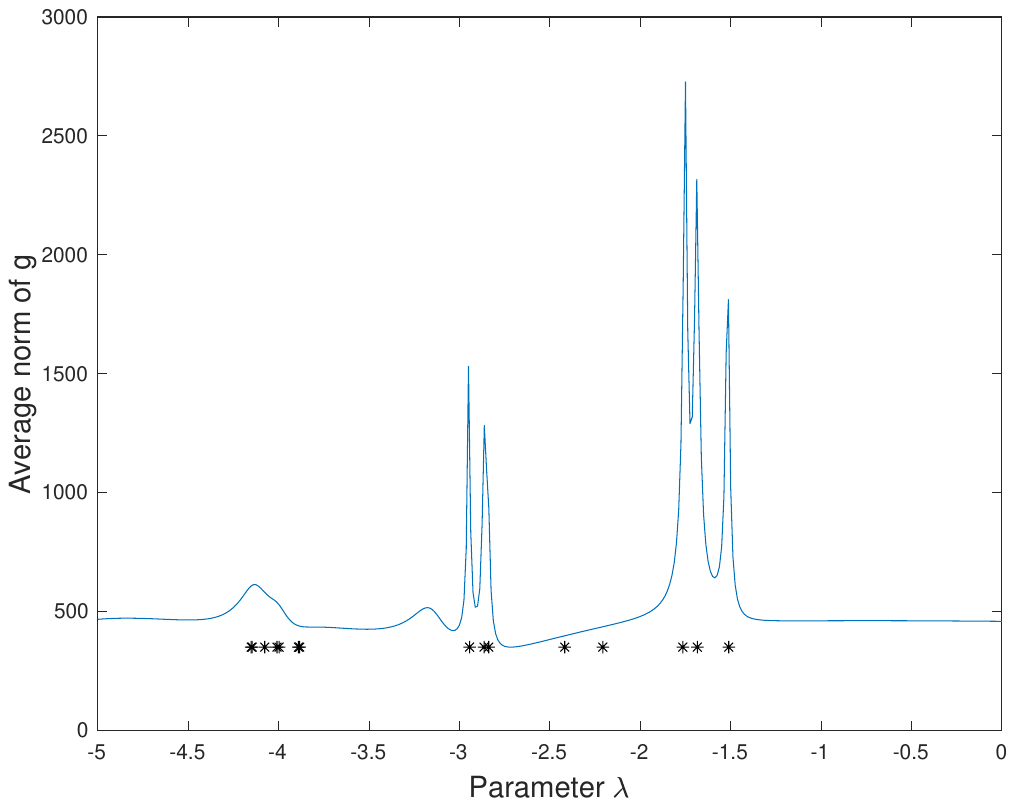}}
\end{center}
\caption{Predicted target signatures and computed $\Sigma$-Steklov eigenvalues for the hemisphere at $\kappa=1.9$.  Left: scalar $\Sigma=0.5iI$. Right: anisotropic $\Sigma$ with $\sigma_1=0.5$ and $\sigma_3=0.4$. In each panel the curve shows the average of the norm of ${\bf g}$ as the parameter $\lambda$ varies, and the $*$ mark eigenvalues computed by FEM.}
\label{hemi}
\end{figure}

\paragraph*{A closed screen:}
A closed spherical screen is a useful test case to check all steps of the algorithm since all problems can be solved analytically using special function expansions.  In the results presented here we assume $\Sigma=\partial B_1$.  Because of constraints on the finite element solver, we choose a modest value $\kappa=1.9$.  We choose $\Sigma$ 
to be the diagonal matrix $\Sigma=(0.5i)I$ resulting in real $\Sigma$-Steklov eigenvalues.  Then we solve the forward problem to generate scattering data which is corrupted by uniformly distributed random noise at each data point introducing 0.15\% error in the computed far field pattern in the relative spectral norm (see \cite{CBMS-kot} for more details).  We also solve the auxiliary problem for 501 choices of $\eta\in [-0.5,1]$.  Results are shown in Fig.~\ref{sphere}.  We see clear detection of the three $\Sigma$-Steklov eigenvalues in this range that
agree well with eigenvalues computed by the FEM (on the vertical scale used in Fig~\ref{sphere}, the leftmost peak is barely visible).
\paragraph*{A hemispherical screen:}
We next consider a hemispherical screen on the surface of the sphere of radius 1.  We first set the scalar parameter $\Sigma=0.5 i I$ and $\kappa=1.9$.  Solving the forward problem by FEM requires a finer mesh near the screen than is needed in the background media as shown in Fig.~\ref{fwd}.  This substantially increases the time for the forward solve, but of course does not affect the computation of target signatures once far field data for the auxiliary problem is computed.
Using data computed by the FEM and corrupted by noise as for the sphere, the resulting predicted target signatures are shown in  the left panel of Fig~\ref{hemi}.  The $\Sigma$-Steklov eigenvalues are changed compared to Fig.~\ref{sphere}.  The results for the leftmost cluster of signatures are smeared out compared to the two other group of eigenvalues (but the vertical scale does not emphasize this cluster).

Next we consider an anisotropic surface
conductivity on the hemispherical screen and take $\Sigma$ and in order to define the anisotropic $\Sigma$ we first define
\[\tilde\Sigma=\left(\begin{array}{ccc}\sigma_{1,1}i&0&0\\0&0.5i&0\\0&0&\sigma_{3,3}i\end{array}\right)
\]
where $\sigma_{1,1}$ and $\sigma_{3,3}$ will be chosen later.
Then for a tangential vector field ${\bf v}$ we set
\begin{equation}
\Sigma {\bf v}=P_{\Gamma}\tilde\Sigma \bfv\label{Sig_proj}
\end{equation}
where $P_{\Gamma}$ denotes projection on to the tangent plane of the sphere at each point of the hemisphere.
For the example in this section, we set $\sigma_{1,1}=0.5$ and $\sigma_{3,3}=0.4$.  Results are shown in the right panel of Fig.~\ref{hemi}.  Although the eigenvalues are changed, the far field only picks up the change in the rightmost eigenvalue.  None-the-less the anisotropy is detected.

\begin{figure}
\begin{center}
\resizebox{0.49\textwidth}{!}{\includegraphics{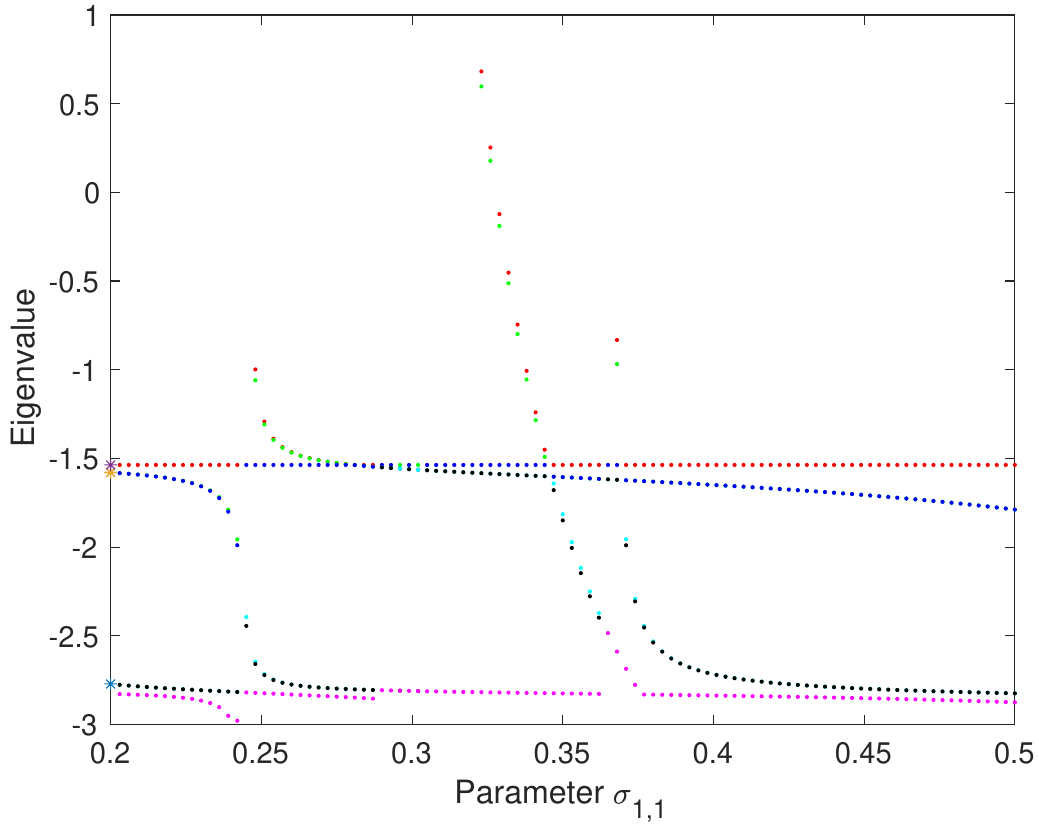}}
\resizebox{0.49\textwidth}{!}{\includegraphics{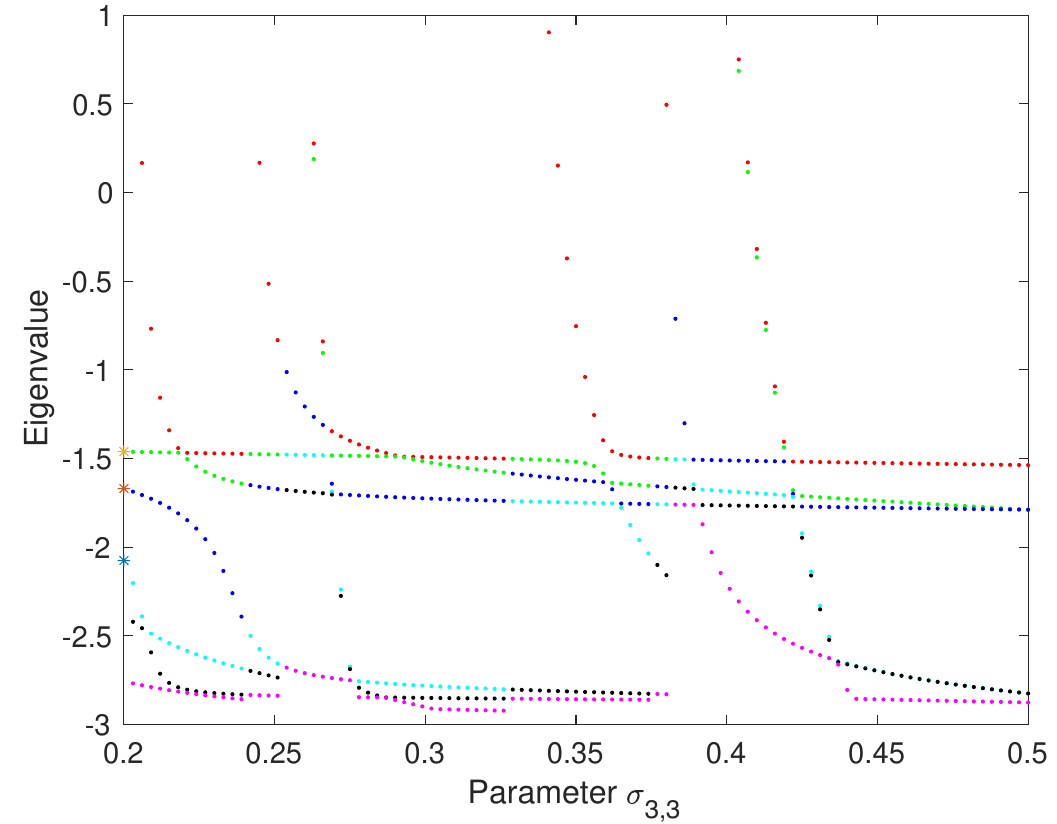}}
\end{center}
\caption{Results of changing parameters in an anisotropic choice of $\Sigma$ for the hemispherical screen. We show changes in the smallest (in magnitude) target signatures as the parameters defining $\Sigma$ given by (\ref{Sig_proj})) vary.  Left panel: we set $\sigma_{3,3}=0.5$ and vary $\sigma_{1,1}$.  Right panel: we set $\sigma_{1,1}=0.5$ and vary $\sigma_{3,3}$.  Eigenvalues for different parameter values are shown as $*$.}\label{follow}
\end{figure}

\paragraph*{Investigating eigenvalues}
The eigensolver can be used to study the effects of changes in $\Sigma$ on the $\Sigma$-Steklov eigenvalues and so predict the sensitivity of the target signature to changes in the surface properties.  
Using the finite element eigensolver discussed in Section~\ref{eigen-num} we can solve the eigenvalue problem for different choices of $\sigma_{1,1}$ and $\sigma_{3,3}$ and follow changes in the target signatures as a function of the surface parameters.
Results are shown in Fig.~\ref{follow}

\section{Conclusion}
We have shown preliminary results for the inverse problem of detecting changes in a thin anisotropic scatterer.  We have provided a general existence theory for the forward problem, as well as a basic uniqueness result for the inverse problem.  We also developed the idea of  $\Sigma$-Steklov eigenvalues as target signatures for the screen.  At present the majority of the theory, and all the numerical results are for purely imaginary surface impedance (a lossless screen).  
Further work is needed to prove the existence of $\Sigma$-Steklov eigenvalues when $\Sigma$ is a complex tensor, and numerical testing in this case is also needed.

\section*{Acknowledgements} 
The research of F.C.  was partially supported  by the  US AFOSR Grant FA9550-23-1-0256 and  NSF Grant DMS 2406313. The research of  P. M. was partially supported by the US AFOSR under grant number FA9550-23-1-0256.  
We thank Professor Timo Lahivaara of the University of Eastern Finland for his help benchmarking the code used in this study.
\bibliographystyle{siamplain}
\bibliography{screen}
\end{document}